\PassOptionsToPackage{unicode}{hyperref}
\PassOptionsToPackage{svgnames}{xcolor}
\documentclass[12pt,a4paper]{article}
\usepackage{authblk}
\usepackage{amsmath}
\usepackage{cases}
\usepackage{natbib}
\usepackage{enumitem}
\usepackage{amsthm}
\usepackage{amssymb}
\usepackage{caption}
\usepackage{graphicx}
\usepackage{algorithm}
\usepackage{algpseudocode}
\newtheorem{theorem}{Theorem}
\newtheorem{assumption}{Assumption}
\newtheorem{definition}{Definition}
\newtheorem{corollary}{Corollary}
\newtheorem{lemma}{Lemma}
\newtheorem{proposition}{Proposition}
\theoremstyle{definition}
\newtheorem{remark}{Remark}

\newcommand{\E}{\mathbb{E}}
\newcommand{\PP}{\mathbb{P}}
\newcommand{\ind}{\mathbb{I}}
\newcommand{\cov}{{\rm Cov}}
\newcommand{\supp}{\mathrm{supp}}
\newcommand{\sign}{\mathrm{sign}}
\newcommand{\bW}{{\mathbf W}}
\newcommand{\bX}{{\mathbf X}}
\newcommand{\bZ}{{\mathbf Z}}
\newcommand{\btheta}{\boldsymbol{\theta}}
\newcommand{\bxi}{\boldsymbol{\xi}}

\usepackage[T1]{fontenc}
\usepackage{lmodern}
\usepackage{microtype}
\UseMicrotypeSet[protrusion]{basicmath}
\usepackage{parskip}
\usepackage{xcolor}
\makeatletter
\def\maxwidth{\ifdim\Gin@nat@width>\linewidth\linewidth\else\Gin@nat@width\fi}
\def\maxheight{\ifdim\Gin@nat@height>\textheight\textheight\else\Gin@nat@height\fi}
\makeatother
\setkeys{Gin}{width=\maxwidth,height=\maxheight,keepaspectratio}
\usepackage[a4paper,margin=1in]{geometry}
\usepackage{setspace}

\usepackage{hyperref}
\usepackage{bookmark}

\hypersetup{
  pdftitle={Learning Sparse Support under Differential Privacy: Adaptive Algorithms and Minimax Limits},
  pdfauthor={Jia Gu; T. Tony Cai},
  pdfkeywords={differential privacy, support recovery, LASSO, minimax limit, propose-test-release},
  colorlinks=true,
  linkcolor={blue},
  filecolor={Blue},
  citecolor={Blue},
  urlcolor={Blue}}

\usepackage{needspace}

\begin{document}
\renewcommand{\baselinestretch}{1}\small\normalsize

\title{\bf Learning Sparse Support under Differential Privacy: Adaptive Algorithms and Minimax Limits}
\author[1]{Jia Gu}
\author[2]{T. Tony Cai\thanks{Corresponding author.}}
\affil[1]{Center for Data Science, Zhejiang University}
\affil[2]{Department of Statistics and Data Science, The Wharton School, University of Pennsylvania}
 \date{}
\maketitle

\begin{abstract}
  We study exact support recovery under $(\epsilon,\delta)$-differential privacy in sparse high-dimensional linear regression. We introduce Saturated Propose-test-release, a
general mechanism that privately releases the output of a discrete selector with
probability one once its stability certificate reaches a finite threshold. Exploiting the
coordinatewise geometry of the LASSO, we construct a computable support-stability
score. The resulting  computationally efficient Saturated LASSO satisfies worst-case $(\epsilon,\delta)$-differential
privacy and achieves exact support recovery with high probability under explicit
regularity and beta-min conditions. Maximizing sparsity-indexed certificates yields
an adaptive procedure requiring no sparsity knowledge and having exactly the
same finite-sample exact-recovery risk as the oracle fixed-sparsity procedure under
common public tuning parameters. We also establish a minimax lower bound
explicitly tracking $\delta$: under its recovery conditions,
Saturated LASSO is minimax optimal up to logarithmic factors in $n$ and $1/\delta$ uniformly over $0 < \delta \leq \epsilon/16$; in
the broad moderate-$\delta$ regime, it further matches the lower-bound $\delta$-dependence.
A complementary information-theoretic construction with known sparsity attains
the lower-bound rates up to constant factors under independent Gaussian design,
at exponential computational cost. Simulations and a semi-synthetic study using
public American Community Survey covariates illustrate the numerical performance
of the proposed procedures.

\end{abstract}


\newpage
\doublespacing

\section{Introduction}

Support recovery, also known as variable selection, is a fundamental problem in statistical and machine learning, particularly in high-dimensional settings where the underlying signal is often sparse. Such problems arise across a wide range of applications, including economics, genomics, and epidemiology. Identifying the relevant predictors can sharpen scientific interpretation, guide subsequent investigation, and improve estimation and prediction. A rich
body of methods with strong theoretical guarantees has been developed for support
recovery in conventional, non-private settings. These works characterize the fundamental
limits of variable selection and propose computationally efficient procedures that achieve
optimal or near-optimal performance under appropriate conditions; see \citet{Fan-JASA-2001, Zhao-JMLR-2006, Zhang-JMLR-2009, Wainwright-TIT-2009, Wainwright-TIT-LASSO-2009, Zhang-AOS-2010, Cai-TIT-2011, Butucea-AoS-2018}.

As statistical learning increasingly relies on data containing sensitive personal information, from health-care records to individual and household profiles, protecting privacy in released statistical results has become essential. Learning algorithms may reveal information about individuals whose records were used for training, including whether a particular individual participated in the dataset \citep{Reza-IEEESP-2017}. Differential privacy (DP), introduced by \citet{Dwork-2006}, provides a rigorous framework for limiting such disclosures while preserving the possibility of useful statistical analysis \citep{Cummings-HDSR-2024}.

Privacy protection is also important in support recovery because the identities of the selected variables are themselves the released statistical output. For example, economic and social research may use confidential survey records to identify factors associated with an outcome of interest. Even if the underlying records are never disclosed, the selected set can depend sensitively on the responses of particular individuals or households and may therefore reveal information about them. The goal of this paper is to develop adaptive and computationally efficient support recovery procedures under $(\epsilon,\delta)$-DP and to characterize their optimality in sparse high-dimensional linear regression.

\subsection{Private sparse support recovery}

Let
$\bZ=(Z_1,\ldots,Z_n)\in\mathcal{Z}^n$ denote a dataset of $n$ i.i.d. observations, where each observation is of the form
$Z=(Y,X^\top)^\top$ and follows the linear regression model
\begin{equation}
Y = X^\top\btheta^* + \varepsilon,
\label{eq: linear regression model}
\end{equation}
where $\btheta^*\in\mathbb{R}^d$ is the unknown population parameter, $X\sim\mathcal{P}_X$ satisfies
$\E(X)=0_d$ and $\cov(X)=\Sigma$, and the noise variable $\varepsilon$ is independent of $X$, has mean zero, and is sub-Gaussian with parameter $\sigma>0$. We assume that the  population parameter $\btheta^* $ is $(s,\theta_{\min})$-sparse, namely $\btheta^* \in   \Theta(s,\theta_{\min})$,
where
\begin{equation}
\Theta(s,\theta_{\min}) = \left\{\btheta\in\mathbb{R}^d: |\supp(\btheta)| = s \mbox{ and } |\theta_j|\geq \theta_{\min} \mbox{ for all } j\in\supp(\btheta)\right\}. \label{eq: sparse class}
\end{equation}
Under the above model, we aim to recover the support $S = \supp(\btheta^*)$ of a population parameter $\btheta^*$ by a selector $\mathcal{M}:\mathcal{Z}^n \to \mathcal{J}\cup\{\perp\}$,  where $\mathcal{J} := \{J: J\subset [d]\}$ and $\perp$ represents the no-release symbol. The selector $\mathcal{M}(\bZ)$ has to satisfy $(\epsilon,\delta)$-differential privacy, which is defined as follows.
\begin{definition}[Differential Privacy]
    A randomized algorithm $\mathcal{M}$ that maps $\mathcal{Z}^n$ to some space $\mathcal{R}$ is $(\epsilon,\delta)$-differentially private (DP) for $\epsilon\geq0$ and $\delta\in[0,1)$ if, for every pair of neighboring datasets $\bZ,\bZ^\prime\in \mathcal{Z}^n$ that differ by one individual datum and every measurable $U\subset \mathcal{R}$,
    \[
\PP(\mathcal{M}(\bZ)\in U ) \leq e^{\epsilon} \cdot \PP(\mathcal{M}(\bZ^\prime)\in U) + \delta,
    \]
    where the probability measure $\PP$ is induced by the randomness of $\mathcal{M}$ only.
\end{definition}
Two datasets $\bZ$ and $\bZ^\prime$ are called neighbors, written
$\bZ\sim\bZ^\prime$, if they differ in exactly one observation. More generally, the Hamming distance
$D_H(\bZ,\bW)$ denotes the minimum number of observation replacements needed to transform
$\bZ$ into $\bW \in \mathcal{Z}^n$. The privacy guarantee considered here is the standard worst-case guarantee: it must hold uniformly over all possible datasets, independently of the data-generating distribution. Hence, the privacy guarantee is model-free. Statistical assumptions enter the analysis of recovery performance, including the sufficient and necessary signal strength conditions developed below.
The performance of $\mathcal{M}(\bZ)$ is measured by the indicator function on the event of exact support recovery error:
$\rho(\mathcal{M}(\bZ),S) = \ind\left(\mathcal{M}(\bZ) \not = S\right)$. This error metric has been used in conventional non-private support recovery \citep{Wainwright-TIT-LASSO-2009}.

\subsection{Main contributions}
\label{sec:main-contributions}

The main contributions of this paper are as follows.

\begin{figure}[!ht]
    \centering
    \includegraphics[width=1.0\linewidth]{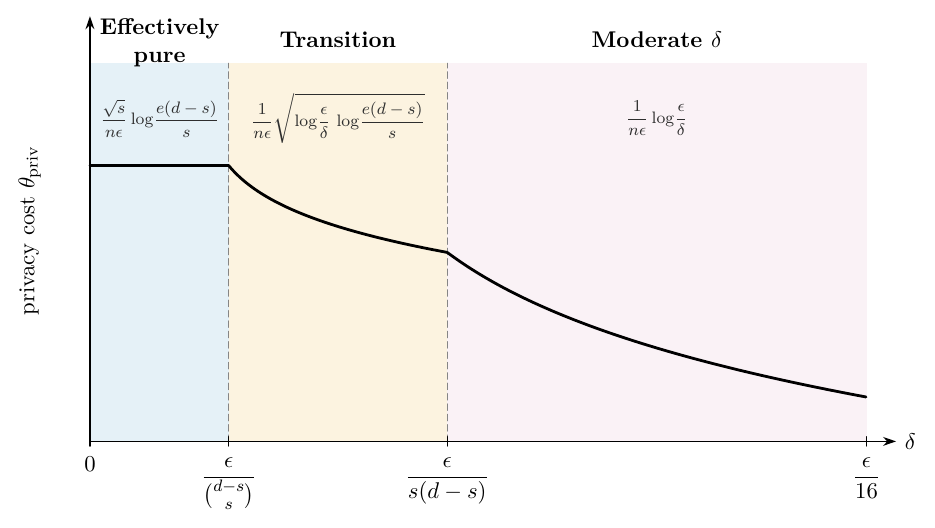}
    \caption{Three regimes of the privacy term in Theorem~\ref{thm:phase-transition}. The maximizing swap scale is $r=1$ in the moderate-$\delta$ regime and $r=s$ in the effectively pure regime, with progressively larger swap scales becoming relevant through the transition regime. The simplified transition formula uses $\log\{e(d-s)/s\}\geq5\log(es)$. }
    \label{fig:theta-priv-phase-transition}
\end{figure}

\begin{enumerate}[label=(\alph*)]
\item \textbf{Saturated PTR and KKT-certified support recovery.}
We develop Saturated
LASSO, which satisfies worst-case $(\epsilon,\delta)$-DP and achieves exact support recovery with
high probability under explicit regularity and minimum signal strength conditions.
The procedure combines a computable stability certificate obtained from coordinatewise perturbation bounds for the LASSO estimate margins with Saturated PTR, a
Propose-test-release (PTR) procedure whose release probability reaches one at a
finite stability threshold. Finite saturation removes the additional stability margin
required by conventional PTR.

\item \textbf{Adaptation to unknown sparsity without an additional  cost.}
We develop an adaptive version of Saturated LASSO that does not require knowledge of the true sparsity level and has exactly the same
finite-sample exact-recovery probability as the corresponding known-sparsity procedure with otherwise identical public tuning parameters. Consequently, adaptation
incurs no additional privacy loss in terms of the minimum signal strength requirement.
This is achieved by aggregating stability certificates for a common LASSO fit and
applying private release once.

\item \textbf{\(\delta\)-dependent minimax rates.}
We establish a minimax lower bound (see Figure~\ref{fig:theta-priv-phase-transition}) for the signal strength required for exact support recovery by comparing supports obtained by exchanging $r$ active coordinates with $r$ inactive coordinates. Under independent Gaussian design and noise, known sparsity, and $\ell_2$-bounded parameter, we provide a complementary exponential-time construction that attains all three lower-bound orders up to constants, establishing the information-theoretic sharpness of the phase transition within that model.

On the common regular model classes defined in Corollary~\ref{cor:minimax-optimality}, and under the stated recovery conditions, both fixed-sparsity and adaptive Saturated LASSO are minimax optimal up to logarithmic factors in $n$ and $1/\delta$, uniformly over $0<\delta\leq\epsilon/16$. A sharper comparison holds in the broad moderate-$\delta$ regime: the sufficient and necessary signal thresholds match up to the $\sqrt{\log n}$ response-clipping factor under sub-Gaussian noise and up to constants under bounded noise.

\end{enumerate}

\subsection{Related work}

Minimax rates under differential privacy have been established for many statistical estimation problems, including mean estimation and linear regression in low- and high-dimensional regimes \citep{Cai-Wang-Zhang-2021}. A corrected arXiv version of that work refines the dependence on \(\delta\) in lower bounds for sparse mean estimation and sparse regression \citep{cai2026cost}, while \citet{Cai-2023} developed general lower-bound tools, including the score-attack method, for a broader class of estimation problems. These results characterize the privacy cost of parameter estimation but do not directly address exact support recovery, whose target is a discrete selected set rather than a continuous parameter.

Private selection has  been studied in model-agnostic settings. \citet{Bafna-COLT-2017} analyzed the cost of private top-\(k\) selection, and \citet{Steinke-FOCS-2017} established tight lower bounds for differentially private selection. Although closely related to support recovery, these results do not exploit the structure of the regression model.

More directly related work concerns private variable selection. The closest existing work is \citet{Guha-COLT-2013}, who introduced perturbation- and subsampling-stability transformations for private support recovery and applied them to the LASSO. Their perturbation-stability construction uses primal--dual analysis to form computable KKT-based tests involving inactive gradients, restricted curvature, and active coefficient bounds, followed by Laplace PTR. These algorithms are not adaptive to the sparsity level. \citet{Lei-JRSSA-2018} proposed differentially private model-selection procedures based on penalized and constrained likelihood and established consistency under conditions resembling those of their non-private counterparts. More recent work has studied private best subset selection using the exponential mechanism, MCMC, or mixed-integer programming \citep{Roy-2024,Prastakos-NIPS-2025} under pure DP.  \citet{Cai-JRSSB-2026,Cai-arxiv-2025} considered private false discovery rate control rather than exact support recovery.

\subsection*{Organization}

The rest of the paper is organized as follows.  Section~\ref{sec:saturated-ptr-discrete-selectors} develops Saturated PTR for stable discrete selectors. Section~\ref{sec:constructing-saturated-lasso} constructs the KKT-certified Saturated LASSO and its sparsity-adaptive version, and provides high-probability exact support recovery guarantees. Section~\ref{sec: cost of privacy} establishes $\delta$-dependent minimax limits and contrasts support recovery with parameter estimation. Section~\ref {sec: numerical results} reports results of simulation studies and a semi-synthetic case study. Section~\ref{sec: discussion} concludes with a discussion.
Technical details are provided in the supplementary material \citep{Gu-DP-supportrecovery}.

\subsection*{Notations}
For a positive integer $d$, write $[d]=\{1,\ldots,d\}$. For a vector $v\in\mathbb R^d$, let $\supp(v)=\{j:v_j\neq0\}$ denote its support, and let $\|v\|_q$ be the usual $\ell_q$ norm. For a set $A\subset[d]$, write $A^c=[d]\backslash A$, $|A|$ for its cardinality, and $v_A$ for the subvector of $v$ indexed by $A$. For a matrix $M\in\mathbb R^{d\times d}$, let $M_{AB}$ denote the submatrix with rows in $A$ and columns in $B$, and let $\|M\|_\infty=\max_i\sum_j |M_{ij}|$. We use $\E$, $\PP$, and $\cov(\cdot)$ to denote expectation, probability, and covariance, respectively.  For two non-negative sequences $a_n$ and $b_n$, write $a_n\lesssim b_n$ if $a_n\leq Cb_n$ for a universal constant $C>0$, and write $a_n\asymp b_n$ if both $a_n\lesssim b_n$ and $b_n\lesssim a_n$ hold. We also follow the convention that $\log(\epsilon/\delta) = \infty $ if $\delta=0$ and $\epsilon > 0 $.

\section{Saturated PTR for stable discrete selectors}\label{sec:saturated-ptr-discrete-selectors}

In Section~\ref{sec: stability score}, we define stable selectors and sub-distance scores, which are two core components of the PTR framework \citep{Dwork-Lei-STOC-2009}. In Section~\ref{sec:saturated-ptr}, we introduce a pointwise maximized
release mechanism that releases the  output of a discrete selector privately with probability one on sufficiently stable datasets.

\subsection{Stable selectors and sub-distances}\label{sec: stability score}
In this section, we adopt the distance-to-instability framework, which has been used for private variable selection \citep{Guha-COLT-2013} and parameter estimation \citep{Shen-arxiv-2026}, and specialize it to discrete support selectors.  Let $\bZ \sim \bZ^\prime$ be a pair of neighboring datasets, and let $\widehat{S}(\bZ)$ be a non-private support estimator. Define
\begin{equation}
\mathrm{sub}_0 (\widehat{S}) = \left\{\bZ \mid \widehat{S}(\bZ^\prime) = \widehat{S}(\bZ) \quad \mbox{ for all } \quad \bZ^\prime \sim \bZ\right\}. \label{eq: sub 0 set}
\end{equation}
Thus, $\mathrm{sub}_0$ consists of datasets for which the support estimator is locally stable under the replacement of a single observation by an arbitrary data point. Next, we define a distance-type score function $\gamma(\bZ)$ to quantify how far a specific dataset $\bZ$ is from the sensitive set $\mathrm{sub}_0^c=\mathcal{Z}^n\setminus\mathrm{sub}_0$. The most direct choice would be the Hamming distance $D_H(\bZ, \mathrm{sub}_0^c)$. However, explicitly constructing $\mathrm{sub}_0$ or computing the exact Hamming distance to its complement is computationally infeasible in practice. Therefore, we use a tractable lower-bound surrogate instead, which we call a sub-distance.
\begin{definition}[Sub-distance]\label{def: sub-distance}
   A non-negative function $\gamma(\cdot) : \mathcal{Z}^n \mapsto [0,\infty)$ is called a sub-distance with respect to a selector $\widehat{S}: \mathcal{Z}^n \mapsto \mathcal{J}\cup\{\perp\}$, if
\begin{equation}
\{\bZ : \gamma(\bZ) > 0\} \subset \mathrm{sub}_0(\widehat{S}),
\qquad \mbox{and} \qquad |\gamma(\bZ) - \gamma(\bZ')| \leq D_H(\bZ,\bZ')
, \label{eq: sub-distance definition}
\end{equation}
where the region $\mathrm{sub}_0(\widehat{S})$ is defined in (\ref{eq: sub 0 set}).
\end{definition}
The first condition ensures that $\gamma(\bZ)>0$ certifies that the dataset $\bZ$ lies in the stable, or safe, region $\mathrm{sub}_0(\widehat{S})$. The second condition states that $\gamma$ is $1$-Lipschitz with respect to Hamming distance, which is critical for constructing private release algorithms.

\subsection{Saturated PTR}\label{sec:saturated-ptr}

Given a dataset $\bZ$ and a sub-distance score $\gamma$, consider a mechanism $\mathcal{M}(\bZ)$ which releases $\widehat S(\bZ)$ with probability $p(\gamma(\bZ))$ and otherwise withholds the answer and outputs $\perp$, where $p:[0,\infty)\to[0,1]$ is nondecreasing.  The intuition is that $\gamma(\bZ)$ measures how deep $\bZ$ lies inside the stable region, so a larger score should allow a higher release probability. We therefore seek to maximize $p(\cdot)$ subject to $(\epsilon,\delta)$-DP. Let $L(\bZ):=\lceil\gamma(\bZ)\rceil$ and write $p(\gamma(\bZ))=q_{L(\bZ)}$ for a nondecreasing sequence $\{q_k\}_{k\geq 0}$. Since $\gamma$ is $1$-Lipschitz, $L$ changes by at most one between neighboring datasets. We can thus optimize $\{q_k\}_{k\geq 0}$ under privacy constraints at adjacent indices $k\leftrightarrow k+1$. For neighboring datasets with different selector outputs, Definition~\ref{def: sub-distance} gives zero scores, so $q_0\leq\delta$ ensures privacy. Taking $q_0=\delta$ and maximizing each $q_{k}$ recursively gives

\begin{equation}
q_{k+1}=\min\left\{1,e^\epsilon q_k+\delta,1-e^{-\epsilon}(1-q_k-\delta)\right\},\qquad k=0,1,2,\ldots.
\label{eq:saturated-ptr-recurrence}
\end{equation}
We  define the mechanism in Algorithm~\ref{alg:saturated-ptr}, and provide its privacy guarantee in  Theorem~\ref{thm:saturated-ptr-privacy}. The pointwise  optimality and finite saturation property of \eqref{eq:saturated-ptr-recurrence} are provided in Proposition~\ref{prop:finite-saturation}.
\begin{algorithm}[!ht]
\caption{Saturated PTR}
\label{alg:saturated-ptr}
\small
\begin{algorithmic}[1]
\Require Deterministic selector $\widehat S:\mathcal Z^n\to\mathcal J\cup\{\perp\}$; sub-distance $\gamma$; privacy parameters $\epsilon>0$ and $\delta\in(0,1)$.
\Ensure Private output $\widehat S_{\rm priv}\in\mathcal J\cup\{\perp\}$.
\State Compute $L(\bZ)\gets\lceil\gamma(\bZ)\rceil$ and the release probability $q_{L(\bZ)}$ from \eqref{eq:saturated-ptr-recurrence}.
\State Release $\widehat S(\bZ)$ with probability $q_{L(\bZ)}$; otherwise release $\perp$.
\end{algorithmic}
\end{algorithm}
\begin{theorem}[Privacy guarantee of Saturated PTR]\label{thm:saturated-ptr-privacy}
Let $\epsilon>0$ and $0<\delta<1$. If $\gamma(\cdot)$ is a sub-distance for $\widehat S$, then Algorithm~\ref{alg:saturated-ptr} is $(\epsilon,\delta)$-DP for all $\bZ \in \mathcal{Z}^n$.
\end{theorem}

\begin{proposition}\label{prop:finite-saturation}
Let $\epsilon>0$, $0<\delta<1$, and let $\{q_k\}_{k\geq0}$ be defined by \eqref{eq:saturated-ptr-recurrence} with $q_0=\delta$. Among all nondecreasing sequences $\{u_k\}_{k\geq0}\subset[0,1]$ with $u_0\leq\delta$ satisfying
\[
u_{k+1}\leq e^\epsilon u_k+\delta,
\qquad
1-u_k\leq e^\epsilon(1-u_{k+1})+\delta,
\qquad k\geq0,
\]
the sequence $\{q_k\}_{k\geq 0}$ is pointwise largest: $u_k\leq q_k$ for every $k\geq0$. Moreover, the saturation index
$K_{\epsilon,\delta}:=\min\{k\geq0:q_k=1\}$ is finite, and $q_k=1$ for every
$k\geq K_{\epsilon,\delta}$. If $0 < \epsilon \leq 1 $ and $0< \delta \leq \epsilon/16$, then
\[
K_{\epsilon,\delta}\leq\frac{4}{\epsilon}\log\!\left(\frac{\epsilon}{\delta}\right).
\]
\end{proposition}

  \begin{remark}
Recurrence~\eqref{eq:saturated-ptr-recurrence} also appears in private binomial testing
\citep{Awan-NeurIPS-2018} and partition selection
\citep{Desfontaines-POPETS-2022}.
Both the count-based statistic and the stability score $\gamma$ change by at most one under data replacement.  Saturated PTR therefore transfers this construction from private testing to the PTR framework: rather than applying a generic randomized test to the stability score, it directly chooses the largest
release probability permitted by the adjacent-layer privacy constraints. To our knowledge, this is the first explicit PTR formulation of this
optimal release rule for general discrete selectors equipped with
sub-distance scores.
\end{remark}

Finite saturation has a direct consequence for exact support recovery: once the certificate reaches $K_{\epsilon,\delta}$, which depends only on the privacy parameters, the selected support is released with probability one. Conventional Laplace PTR \citep{Guha-COLT-2013, Dwork-2014} and Logistic PTR \citep{Shen-arxiv-2026} have a positive withholding probability at every finite score and require an additional $c\log(d)/\epsilon$ stability margin to reduce the mechanism-failure probability to $d^{-c}$; also see Section~\textcolor{Blue}{S2} of the supplementary material. NetPTR \citep{Wang-Shen-arxiv-2026} uses a truncated logistic release rule that also saturates at a finite certificate threshold, but does not follow the pointwise optimal recurrence \eqref{eq:saturated-ptr-recurrence}. Combining Saturated PTR with the KKT certificate developed in the next section removes this additional $\log(d)/\epsilon$ margin from the sufficient beta-min condition for high-probability $(\epsilon,\delta)$-DP exact support recovery.

\begin{remark}[Interpretation of non-release] The output $\perp$ represents abstention and is distinct from the empty
support $\emptyset$. It indicates that the observed stability
certificate is insufficient to guarantee release at the prescribed privacy level, rather than that the true support is empty. Throughout, a no-release output $\perp$ is counted as a failure of exact support recovery.
Operationally, this reject option can be preferable to always releasing
a potentially inaccurate support estimate, as in noisy selection. We will discuss the implications of the non-release further in the case study presented in Section \ref{sec: real-data}.
\end{remark}

\section{Saturated LASSO: construction, adaptation, and recovery guarantees}\label{sec:constructing-saturated-lasso}

We now develop the Saturated LASSO algorithm for exact support recovery. Section \ref{sec: perturbation and stability} briefly reviews a primal-dual characterization of the LASSO support and identifies the support-changing KKT boundaries and controls their movement under data replacement. Section~\ref{sec:saturated-lasso-stability-score} then translates the bounds into a smooth, computable stability certificate for Algorithm \ref{alg:saturated-ptr}, yielding the fixed-sparsity worst-case $(\epsilon,\delta)$-DP Saturated LASSO.
Section~\ref{sec: adaptive support recovery} achieves adaptivity to the sparsity level $s$ without privacy composition or statistical loss. Section~\ref{sec:upper-guarantees} establishes high-probability exact support recovery under explicit model-dependent regularity assumptions and sufficient beta-min conditions.

\subsection{KKT geometry and perturbation stability}\label{sec: perturbation and stability}
We first briefly review the LASSO algorithm and its induced support estimate, which provides the main geometric intuition underlying our $(\epsilon,\delta)$-DP support recovery algorithm. We start with data clipping. Let $B_X, B_Y > 0 $ be public clipping radii, and define $\widetilde X_{ij}=Q_{B_X}(X_{ij})=\max\{-B_X,\min(X_{ij},B_X)\}$ and $\widetilde Y_i=Q_{B_Y}(Y_i)$.
Let $\widetilde X_i=(\widetilde X_{i1},\ldots,\widetilde X_{id})^\top$ and define $\widehat\Sigma(\bZ):=\frac{1}{n}\sum_{i=1}^n \widetilde X_i\widetilde X_i^\top$ and
$\widehat a(\bZ):=\frac{1}{n}\sum_{i=1}^n \widetilde X_i\widetilde Y_i$. For neighboring datasets $\bZ\sim\bZ^\prime$,
\begin{equation}
\|\widehat a(\bZ)-\widehat a(\bZ^\prime)\|_\infty
\le \eta_a:=\frac{2B_XB_Y}{n},
\qquad
\max_{j,k}|\widehat\Sigma_{jk}(\bZ)-\widehat\Sigma_{jk}(\bZ^\prime)|
\le \eta_\Sigma:=\frac{2B_X^2}{n}.
\label{eq:eta}
\end{equation}
The LASSO solves the following $\ell_1$-penalized optimization problem
\begin{equation}
\widehat\btheta^\lambda(\bZ)\in\arg\min_{b\in\mathbb{R}^d}
\left\{\frac{1}{2}b^\top\widehat\Sigma(\bZ)b-\widehat a(\bZ)^\top b+\lambda\|b\|_1\right\}
\label{eq:lasso}.
\end{equation}
Throughout, $\widehat\btheta^\lambda(\bZ)$ denotes the LASSO minimizer with smallest Euclidean norm. This fixed deterministic measurable selection rule is used for all candidate sparsity levels.
Then, a natural support estimator is $\widehat S_\lambda(\bZ):=\supp(\widehat\btheta^\lambda(\bZ))$. For any $A \subset [d]$, Lemma \ref{lem: PDW} provides a set of conditions ensuring $\widehat{S}_{\lambda}(\bZ) = A $. When  $A$ is set as $S$, this reduces to the  primal-dual witness characterization of correct support recovery  \citep{Wainwright-TIT-LASSO-2009}.

\begin{lemma}\label{lem: PDW}
Fix $A\subset[d]$ and $q_A\in\{\pm1\}^{|A|}$, and assume $\widehat\Sigma_{AA}(\bZ)$ is invertible. Define the active oracle
\begin{equation}
\vartheta_A(\bZ;q_A):=\widehat\Sigma_{AA}(\bZ)^{-1}\{\widehat a_A(\bZ)-\lambda q_A\}.
\label{eq:active-oracle}
\end{equation}
Further define the active and inactive KKT margins
\[
m_A^{\rm primal}(\bZ):=\min_{j\in A}q_{A,j}\vartheta_{A,j}(\bZ),
\qquad
m_A^{\rm dual}(\bZ):=\lambda-\max_{j\notin A}|\widehat a_j(\bZ)-\widehat\Sigma_{jA}(\bZ)\vartheta_A(\bZ)|.
\]
If both margins are positive, the vector $\beta\in\mathbb{R}^d$ defined by $\beta_A=\vartheta_A(\bZ)$ and $\beta_{A^c}=\mathbf0$ is the unique minimizer of \eqref{eq:lasso}. In particular, \(\widehat S_\lambda(\bZ)=A\) and \(\sign(\widehat\btheta_A^\lambda(\bZ))=q_A\).
\end{lemma}

 The primal margin $m_A^{\rm primal}$ measures the distance to an active coefficient crossing zero, whereas the dual margin $m_A^{\rm dual}$ measures the distance to an inactive coordinate violating the LASSO dual-feasibility condition. Hence, away from the primal and dual KKT boundaries, the estimated support is locally constant. Establishing stability guarantees therefore requires controlling how rapidly these boundaries can move under data replacements.

To this end, we introduce three empirical regularity buffers for the cross-covariance, active oracle, and inverse stability, respectively. For $A \subset [d]$, define
\[
\Psi_A^{\rm Cross}(\bZ) := \|\widehat{\Sigma}_{A^c,A}\|_\infty = \max_{j\not\in A}\|\widehat{\Sigma}_{jA}(\bZ)\|_1, \quad \Psi_A^{\rm Param}(\bZ) := \|\vartheta_A(\bZ)\|_1, \quad \mbox{ and }
\]
\[
\Psi_A^{\rm Inv}(\bZ) := \begin{cases}
    \|\widehat{\Sigma}_{AA}(\bZ)^{-1}\|_\infty^{-1}, \quad \mbox{if} \ \widehat{\Sigma}_{AA}(\bZ) \ \mbox{ is invertible}\\
    0, \quad \mbox{ otherwise}
\end{cases},
\]
where $\vartheta_A(\bZ)$ is defined in (\ref{eq:active-oracle}).  The algorithm  takes three public constants $\psi^{\rm Cross} > 0$,  $\psi^{\rm Param} > 0$, $\psi^{\rm Inv} > 0$, and we will provide a selection rule in Corollary~\ref{cor:public-buffer-aggregation} of Section~\ref{sec:upper-guarantees}.
Define
\begin{equation*}
\tau_{\rm primal} :=\psi^{\rm Inv}(\eta_a + \eta_{\Sigma}\psi^{\rm Param}),\quad
\tau_{\rm dual} := \eta_a + \eta_{\Sigma}\psi^{\rm Param} + \psi^{\rm Cross} \tau_{\rm primal},
\end{equation*}
where $\eta_a$ and $\eta_{\Sigma}$ are defined in \eqref{eq:eta}. The quantities $\tau_{\rm primal}$ and $\tau_{\rm dual}$ bound the one-step perturbations of the primal and dual KKT margins, respectively. Theorem~\ref{thm: radius stability} uses these scales to establish support stability over an entire Hamming ball.

\begin{theorem}[Hamming-radius support stability]\label{thm: radius stability}
For integers $r>0$ and $s>0$, let $A=\widehat S_\lambda(\bZ)$. Suppose: (i) $|A|\leq s$; (ii) $\Psi_A^{\rm Inv}(\bZ)>(\psi^{\rm Inv})^{-1}+rs\eta_\Sigma$; (iii) $\Psi_A^{\rm Param}(\bZ)<\psi^{\rm Param}$; (iv) $\Psi_A^{\rm Cross}(\bZ)<\psi^{\rm Cross}-rs\eta_\Sigma$; (v) $m_A^{\rm primal}(\bZ)>r\tau_{\rm primal}$; and (vi) $m_A^{\rm dual}(\bZ)>r\tau_{\rm dual}$. Then, for every $\bW$ satisfying $D_H(\bW,\bZ)\leq r$, \(\widehat S_\lambda(\bW)=\widehat S_\lambda(\bZ)\).
\end{theorem}

\subsection{A smooth stability score and Saturated LASSO}\label{sec:saturated-lasso-stability-score}
The perturbation bounds above certify support stability over a fixed Hamming radius, but a sub-distance score (Definition~\ref{def: sub-distance}) must itself be $1$-Lipschitz. The unbuffered constraint $\Psi_A^{\rm Param}<\psi^{\rm Param}$ can cross its threshold after a single observation replacement, causing the score to drop discontinuously to zero. Accordingly, we protect $\Psi_A^{\rm Param}(\bZ)<\psi^{\rm Param}$ with an explicit buffer. For a selected support $A=\widehat S_\lambda(\bZ)$ with $1\leq |A|\leq s$ and invertible $\widehat\Sigma_{AA}(\bZ)$, define
\begin{eqnarray*}
\Delta_A(\bZ) &:=& \min\Bigg\{\frac{\Psi_A^{\rm Inv}(\bZ) - (\psi^{\rm Inv})^{-1}}{s\eta_{\Sigma}},\frac{\psi^{\rm Param} - \Psi_A^{\rm Param}(\bZ)}{s\tau_{\rm primal}},\\
&& \hspace{0.45in}\frac{\psi^{\rm Cross} - \Psi_A^{\rm Cross}(\bZ)}{s\eta_{\Sigma}}, \frac{m_A^{\rm primal}(\bZ)}{\tau_{\rm primal}}, \frac{m_A^{\rm dual}(\bZ)}{\tau_{\rm dual}}\Bigg\}
\end{eqnarray*}
We then define the following score on every dataset
\begin{equation}
\gamma(\bZ) := \begin{cases}
0 \quad \mbox{if} \ |A| > s \ \mbox{or} \ A =\emptyset, \ \mbox{or} \ \widehat{\Sigma}_{AA} \ \mbox{is singular} \\
    \max\left\{\lceil \Delta_A(\bZ) \rceil-1, 0\right\} \quad \mbox{otherwise}
\end{cases} . \label{eq: sub-distance}
\end{equation}
Whenever $\Delta_A(\bZ)>r$ for a positive integer $r$, we have $\gamma(\bZ)\geq r$. Each term in the minimum represents a normalized margin to one of the stability conditions. In particular, a positive score implies positive primal and dual KKT margins and hence uniqueness of the LASSO minimizer by Lemma~\ref{lem: PDW}; nonunique cases therefore receive score zero.
\begin{theorem}[Sub-distance property] \label{thm: sub-distance property}The statistic $\gamma(\cdot)$ defined in (\ref{eq: sub-distance}) is a valid sub-distance measure with respect to the support selector $\widehat{S}_{\lambda}(\cdot)$.
\end{theorem}

The construction and analysis of \eqref{eq: sub-distance} build on the KKT-based stability tests of \citet{Guha-COLT-2013}. Their inactive-gradient test uses the threshold $\lambda/2$ in our normalization, whereas $m_A^{\rm dual}$ uses the full threshold $\lambda$. We bound active-oracle perturbations directly in the $\ell_\infty$ norm through the empirical inverse, rather than through a restricted-curvature $\ell_2$ bound. For fixed public buffers, the resulting margin perturbation scales $\tau_{\rm primal}$ and $\tau_{\rm dual}$ have no explicit sparsity factor. Explicit inverse, cross-covariance, and parameter-norm buffers control the empirical regularity conditions under data replacement and yield a globally $1$-Lipschitz support-stability certificate.

Applying Saturated PTR to the selector $\widehat S_\lambda$ with this score gives Saturated LASSO, presented in Algorithm~\ref{alg:saturated-lasso}. At score zero, the procedure releases the selected support with the baseline probability $q_0=\delta$ and returns $\perp$ with probability $1-\delta$. The privacy guarantee in Theorem~\ref{thm:saturated-lasso-privacy} therefore follows from Theorems~\ref{thm:saturated-ptr-privacy} and~\ref{thm: sub-distance property}.
\begin{algorithm}[!ht]
\caption{Saturated LASSO}
\label{alg:saturated-lasso}
\small
\begin{algorithmic}[1]
\Require Dataset $\bZ=(Z_1,\ldots,Z_n)$ with $Z_i=(Y_i,X_i^\top)^\top$;
privacy parameters $\epsilon>0$ and $\delta\in(0,1)$;
clipping radii $B_X,B_Y>0$;
LASSO tuning parameter $\lambda>0$; sparsity level $s$;
public constants $\psi^{\rm Param}>0$, $\psi^{\rm Inv}>0$, and $\psi^{\rm Cross}>0$.
\Ensure Private output $\widehat S_{\rm priv}\in 2^{[d]}\cup\{\perp\}$.

    \State Clip the covariates and response $\widetilde X_{ij}
    \gets
    \operatorname{sign}(X_{ij})
    \min\{|X_{ij}|,B_X\}$ and $\widetilde Y_i
    \gets
    \operatorname{sign}(Y_i)
    \min\{|Y_i|,B_Y\}$ for $ 1 \leq i \le n$ and $1\leq j \leq d$.

\State Compute the clipped empirical covariance and cross-moment: $\widehat\Sigma
\gets
\frac{1}{n}\sum_{i=1}^n \widetilde X_i\widetilde X_i^\top$ and $\widehat a
\gets
\frac{1}{n}\sum_{i=1}^n \widetilde X_i\widetilde Y_i$.
\State Solve the LASSO problem in~\eqref{eq:lasso} via the selection rule specified after~\eqref{eq:lasso} to obtain
$\widehat\btheta^\lambda(\bZ)$.
\State Set $A\gets \widehat S_\lambda(\bZ)$ and $q_A\gets \operatorname{sign}(\widehat\btheta^\lambda_A(\bZ))$.
\State Compute the neighbor perturbation scales: $\eta_a\gets \frac{2B_XB_Y}{n}$ and $\eta_\Sigma\gets \frac{2B_X^2}{n}$.
\State Compute the perturbation scales: $\tau_{\rm primal}
\gets
\psi^{\rm Inv}\bigl(\eta_a+\eta_\Sigma \psi^{\rm Param}\bigr)$ and  $\tau_{\rm dual}
\gets
\eta_a+\eta_\Sigma \psi^{\rm Param}+\psi^{\rm Cross}\tau_{\rm primal}$.

\State Compute the stability sub-distance score $\gamma(\bZ)$ defined in (\ref{eq: sub-distance}).

\State Use Saturated PTR defined in Algorithm~\ref{alg:saturated-ptr} to release an estimate of the support privately.
\end{algorithmic}
\end{algorithm}

\begin{theorem}\label{thm:saturated-lasso-privacy}
Algorithm \ref{alg:saturated-lasso} is $(\epsilon,\delta)$-DP for all $\bZ \in \mathcal{Z}^n$.
\end{theorem}
The guarantee is worst-case: it follows only from the valid KKT sub-distance and Saturated PTR, and therefore holds for every dataset without requiring the regression model assumption (\ref{eq: linear regression model}) or any random design assumptions imposed in Section \ref{sec:upper-guarantees}.

\subsection{Cost-free sparsity adaptation}\label{sec: adaptive support recovery}

Saturated LASSO (Algorithm \ref{alg:saturated-lasso}) uses the sparsity level only in its stability certificate. For any fixed public penalty level $\lambda$ and a fixed deterministic rule for selecting a LASSO minimizer, the support selector $\widehat S_\lambda$ is common to all candidate sparsity levels. This structure permits the candidate certificates to be aggregated without changing the selected support. More generally, let $\widehat S:\mathcal Z^n\mapsto\mathcal J\cup\{\perp\}$ be a deterministic selector and let $\{\gamma_t(\bZ):t\in\mathcal T\}$ be alternative sub-distance certificates for this same statistical decision. Each positive $\gamma_t(\bZ)$ certifies that $\widehat S(\bZ)$ remains unchanged in a Hamming neighborhood of $\bZ$, and a larger score certifies a larger neighborhood. It is therefore natural to keep the strongest available certificate by taking the pointwise maximum. The next proposition shows that this maximum preserves both defining properties of a sub-distance.

\begin{proposition}[Maximum of sub-distances]\label{prop:maximum-of-sub-distances}
Let $\mathcal T$ be a finite nonempty candidate set and let $\widehat S:\mathcal Z^n\to\mathcal J\cup\{\perp\}$ be a deterministic selector. For every $t\in\mathcal T$, suppose that $\gamma_t:\mathcal Z^n\to[0,\infty)$ is a sub-distance for the same selector $\widehat S$ satisfying \eqref{eq: sub-distance definition}. Define
\begin{equation}
\gamma_{\rm ad}(\bZ):=\max_{t\in\mathcal T}\gamma_t(\bZ).
\label{eq: adaptive maximum score}
\end{equation}
Then $\gamma_{\rm ad}$ is a sub-distance for $\widehat S$.
\end{proposition}
Because the candidate scores are sub-distances for the same selector, Proposition \ref{prop:maximum-of-sub-distances} ensures that their maximum is itself a valid sub-distance. A single application of Saturated PTR with this maximum score therefore gives the privacy guarantee without composition cost. Let $\mathcal{T}=\{1,\ldots,n\}$
be the candidate set, and let $\lambda$ be a public LASSO penalty parameter.
All quantities below are computed after applying the same clipping step as in Algorithm~\ref{alg:saturated-lasso}. Define $A = \widehat{S}_{\lambda}(\bZ)$ and, when $A\neq\emptyset$, set $q_A=\operatorname{sign}(\widehat\btheta_A^\lambda(\bZ))$. For each candidate $t$, let $\mathcal B_t(\bZ)$ denote the event that $|A|>t$, or $A=\emptyset$, or $\widehat\Sigma_{AA}(\bZ)$ is singular. On $\mathcal B_t(\bZ)^c$, define
\[
\begin{aligned}
\Delta_{A,t}(\bZ)
&:=
\min\biggl\{
\frac{\Psi_{A}^{\rm Inv}(\bZ)-(\psi^{\rm Inv})^{-1}}{t\eta_{\Sigma}},
\frac{\psi^{\rm Param}-\Psi_{A}^{\rm Param}(\bZ)}{t\tau_{{\rm primal}}},
\\
&\hspace{0.45in}
\frac{\psi^{\rm Cross}-\Psi_{A}^{\rm Cross}(\bZ)}{t\eta_{\Sigma}},
\frac{m_{A}^{\rm primal}(\bZ)}{\tau_{{\rm primal}}},
\frac{m_{A}^{\rm dual}(\bZ)}{\tau_{{\rm dual}}}
\biggr\}.
\end{aligned}
\]
Then, define
\begin{equation}
\gamma_t(\bZ)
:=
\begin{cases}
0, \quad \mbox{if } \mathcal B_t(\bZ) \mbox{ holds},\\
\max\{\lceil\Delta_{A,t}(\bZ)\rceil-1,0\},\quad \mbox{otherwise}.
\end{cases}
\label{eq: adaptive candidate score}
\end{equation}
By Theorem \ref{thm: sub-distance property}, for each fixed $t$, the score $\gamma_t$ is a valid sub-distance for the common selector $\widehat{S}_{\lambda}$. Therefore, Proposition \ref{prop:maximum-of-sub-distances} shows that $\gamma_{\rm ad}$ in \eqref{eq: adaptive maximum score} is also a valid sub-distance for $\widehat S_\lambda$. The next proposition ensures that the adaptive score admits a simple closed-form expression.
\begin{proposition}\label{prop: adaptive score closed form}
Let $A=\widehat S_\lambda(\bZ)$ and $a=|A|$. If $A=\emptyset$, $a>n$, or $\widehat\Sigma_{AA}(\bZ)$ is singular, then $\gamma_{\rm ad}(\bZ)=0$. Otherwise, $\gamma_{\rm ad}(\bZ)=\gamma_a(\bZ)$.
\end{proposition}
Thus, although $\gamma_{\rm ad}$ is the maximum of a public family of fixed-$t$ sub-distances, it can be computed by substituting the observed support size $a=|\widehat S_\lambda(\bZ)|$ into the candidate score when $1\leq a\leq n$ and $\widehat\Sigma_{AA}(\bZ)$ is nonsingular. Its global smoothness follows from the maximum representation in \eqref{eq: adaptive maximum score}, rather than from treating the data-dependent value $a$ as fixed.
Algorithm \ref{alg:adaptive-saturated-lasso} formally defines the  Adaptive Saturated LASSO procedure, whose theoretical guarantees are provided in Theorem \ref{thm:adaptive-oracle-equivalence}.

\begin{algorithm}[!ht]
\caption{Adaptive Saturated LASSO}
\label{alg:adaptive-saturated-lasso}
\small
\begin{algorithmic}[1]
\Require Dataset $\bZ=(Z_1,\ldots,Z_n)$ with $Z_i=(Y_i,X_i^\top)^\top$;
privacy parameters $\epsilon>0$ and $\delta\in(0,1)$;
clipping radii $B_X,B_Y>0$;
candidate set $\mathcal T=\{1,\ldots,n\}$;
penalty level $\lambda$;
public constants $\psi^{\rm Param}>0$, $\psi^{\rm Inv}>0$, and $\psi^{\rm Cross}>0$.
\Ensure Private adaptive output $\widehat S_{\rm priv}^{\rm ad}\in 2^{[d]}\cup\{\perp\}$.

\State Clip the covariates and response, and  compute $\widehat\Sigma$, $\widehat a$, $\eta_a$, and $\eta_\Sigma$ as in Algorithm \ref{alg:saturated-lasso}.
\State Solve the LASSO problem in~\eqref{eq:lasso} with tuning parameter $\lambda$ via  the selection rule specified after~\eqref{eq:lasso} to obtain $\widehat\btheta^{\lambda}(\bZ)$; set $A\gets\widehat S_\lambda(\bZ)$ and, if $A\neq\emptyset$, $q_A\gets\operatorname{sign}(\widehat\btheta_A^\lambda(\bZ))$.
\State Compute $\tau_{\rm primal}$ and $\tau_{\rm dual}$ as in Algorithm \ref{alg:saturated-lasso}.
\State Compute $\gamma_{\rm ad}(\bZ)=\max_{t\in\mathcal T}\gamma_t(\bZ)$ using \eqref{eq: adaptive candidate score}; equivalently, use Proposition \ref{prop: adaptive score closed form}.
\State Use Saturated PTR in Algorithm~\ref{alg:saturated-ptr} with selector $\widehat S_\lambda$ and score $\gamma_{\rm ad}$ to output $\widehat S_{\rm priv}^{\rm ad}$.
\end{algorithmic}
\end{algorithm}

\begin{theorem}
\label{thm:adaptive-oracle-equivalence}
Algorithm \ref{alg:adaptive-saturated-lasso} is $(\epsilon,\delta)$-DP. Let $S=\supp(\btheta^*)$ and $s=|S| \in \mathcal{T}$.  Then,
\begin{equation}
\PP\left(\widehat S_{\rm priv}^{\rm ad}=S\mid\bZ\right)
=
\PP\left(\widehat S_{\rm priv}=S\mid\bZ\right)
=
\ind\{\widehat S_\lambda(\bZ)=S\}q_{\lceil\gamma_s(\bZ)\rceil},
\label{eq:adaptive-pointwise-oracle-equivalence}
\end{equation}
where $\widehat S_{\rm priv}^{\rm ad}$ denotes the adaptive selector in Algorithm \ref{alg:adaptive-saturated-lasso}, $\widehat S_{\rm priv}$ denotes the non-adaptive private selector in Algorithm \ref{alg:saturated-lasso} with known sparsity level $s$, and $\widehat S_\lambda(\bZ)$ is the corresponding non-private LASSO selector.   Consequently, under the assumptions of Proposition \ref{thm: sufficient beta-min}, the adaptive estimator $\widehat S_{\rm priv}^{\rm ad}$ inherits the same statistical guarantees established therein.
\end{theorem}
 By \eqref{eq:adaptive-pointwise-oracle-equivalence}, the adaptive and oracle fixed-$s$ procedures have identical conditional output distributions when the non-private LASSO recovers $S$; otherwise, neither procedure recovers $S$. Consequently, the two procedures have identical finite-sample exact-recovery probabilities. Therefore, neither a union bound over $\mathcal T$ nor privacy composition is needed and the procedure incurs no cost for adaptation. By contrast,
existing adaptation rules usually run the corresponding algorithm over a
 geometric grid of candidate sparsity levels and compose privacy across the grid, which would  introduce an additional $O(\log n)$ factor in the privacy cost  \citep{Cai-arxiv-2025} under basic composition with equal budget allocation.
\subsection{Exact-recovery statistical guarantees}\label{sec:upper-guarantees}
Having constructed the fixed-sparsity and adaptive procedures and provided the model-free, worst-case privacy guarantee, we now study their exact-recovery performance under the linear regression model \eqref{eq: linear regression model}. By the oracle-equivalence result in Theorem~\ref{thm:adaptive-oracle-equivalence}, it suffices to establish the statistical guarantee for the fixed-sparsity procedure. We first establish a bridge theorem (Theorem~\ref{thm: deterministic beta-min}), which gives data-dependent conditions that guarantee correct exact support recovery. Proposition~\ref{thm: sufficient beta-min} then shows that these
conditions hold with high probability under the statistical model and explicit regularity conditions.
Theorem~\ref{cor: bounded noise rate} then gives the sufficient
beta-min rates under sub-Gaussian and bounded noises.

\begin{theorem}\label{thm: deterministic beta-min}
Let $S$ be the true support of $\btheta^*$, $s = |S|$ and $q_S = \sign(\btheta_S^*)$. Suppose there exists a positive integer $R$ such that
  \[
  \min\left\{ \frac{\Psi_S^{\rm Inv}(\bZ) - (\psi^{\rm Inv})^{-1}}{s\eta_{\Sigma}}, \frac{\psi^{\rm Cross} - \Psi_S^{\rm Cross}(\bZ)}{s\eta_{\Sigma}}, \frac{\psi^{\rm Param} - \Psi_S^{\rm Param}(\bZ)}{s\tau_{\rm primal}}\right\} > R.
  \]
  Also suppose that there exists a constant $\kappa \in (0,1)$ such that $\|\widehat{\Sigma}_{S^c S}\widehat{\Sigma}_{SS}^{-1}\|_\infty \leq 1 -\kappa$. Define $\bxi = \widetilde{\bX}^\top (\widetilde{{{\mathbf Y}}} - \widetilde{\bX}\btheta^*) / n$. If
\[
\lambda > \frac{(2-\kappa)\|\bxi\|_\infty + R \tau_{\rm dual}}{\kappa} \quad \mbox{ and } \quad \theta_{\min} > \psi^{\rm Inv} (\|\bxi\|_\infty+\lambda) + R\tau_{\rm primal},
\]
then we have $\widehat{S}_\lambda(\bZ) = S$ and $\gamma(\bZ) \geq R$, where $\gamma(\cdot)$ is defined in (\ref{eq: sub-distance}).
\end{theorem}
The conditions collect empirical inverse regularity, parameter scale and cross-covariance buffers, irrepresentability, and noise control.
Taking $R=K_{\epsilon,\delta}$
guarantees that the private output is $S$ with probability one conditional on the event.

\begin{assumption}[Regularity conditions]\label{assumption: regularity}
\leavevmode
\begin{enumerate}[label=(\roman*),topsep=0pt]
\item {\rm \bf Bounded parameter.}  The true parameter satisfies
$\|\btheta^*\|_1<a_0$ for some constant $0<a_0<\infty$, and
$\|\btheta^*\|_0=s$.
\item {\rm \bf Bounded design.} There is a constant $B_x<\infty$ such that $\|X\|_\infty<B_x$ almost surely, and $X\sim\mathcal P_X$ is sub-Gaussian with parameter $\sigma_x$.
\item {\rm \bf Incoherence.} The covariance matrix $\Sigma=\cov(X)$ satisfies $\|\Sigma_{SS}^{-1}\|_\infty\leq a_1$, $\sup_j\|\Sigma_{jS}\|_1\leq a_2$, and $\|\Sigma_{S^cS}\Sigma_{SS}^{-1}\|_\infty\leq1-a_3$ for constants $a_1,a_2>0$ and $0<a_3<1$.
\end{enumerate}
\end{assumption}
\begin{remark}
   The  $\ell_1$-boundedness assumption on $\btheta^*$  has also been used for private model selection in low- and high-dimensional regimes \citep{Lei-JRSSA-2018,Roy-2024}. Bounded covariates are standard in private high-dimensional estimation and inference
\citep{Cai-Wang-Zhang-2021}. The incoherence condition holds for particular dependence families after their induced-$\ell_\infty$ inverse and irrepresentability constants are verified; geometric Toeplitz correlation with $|\rho|$ bounded away from one is one such standard example.
\end{remark}

\begin{proposition}[Finite-sample recovery]\label{thm: sufficient beta-min}
Fix $\alpha\in(0,1/4)$, $0<\epsilon\leq1$, and
$0<\delta\leq\epsilon/16$. Under Assumption~\ref{assumption: regularity},
take the model and regularity bounds as fixed public constants.
Fix a public constant $B_X\geq B_x$ and set
$B_Y=\tau_0\sqrt{2\log(12n/\alpha)}$, where
$\tau_0^2=\sigma_x^2a_0^2+\sigma^2$, and
$\psi^{\rm Param}=4a_0$, $\psi^{\rm Inv}=4a_1$, and
$\psi^{\rm Cross}=3a_2$.  Then there exist positive constants $c_0,C_0,C_1$,
depending only on $(a_0,a_1,a_2,a_3,B_x,\sigma_x,\sigma)$, such that, whenever
\[
n\geq
\max\left\{
\frac{s^2\log(12ds/\alpha)}{c_0^2},
\frac{(B_XB_Y+B_X^2)sK_{\epsilon,\delta}}{c_0}
\right\} \quad \mbox{ and }
\]
\[
\lambda
=
C_0\left\{
\sqrt{\frac{\log(12d/\alpha)}{n}}
+
\frac{(B_XB_Y+B_X^2)K_{\epsilon,\delta}}{n}
\right\},
\]
both the fixed-sparsity
and adaptive Saturated LASSO procedures
correctly recover $S$  with probability at least $1-\alpha$ as long as
 $\theta_{\min}\geq C_1\lambda$, where the probability includes the
sample and mechanism randomness. The constants are uniform over the
distributions satisfying the stated bounds and do not depend on
$(n,d,s,\alpha,\epsilon,\delta)$ .
\end{proposition}

\begin{theorem}[Recovery rates]\label{cor: bounded noise rate}
Under Assumption~\ref{assumption: regularity}, use the privacy and
sample-size conditions of Proposition~\ref{thm: sufficient beta-min} and
its public buffers and exact penalty. For fixed $c',c''>0$, take
$\alpha=d^{-c'}+n^{-c''}<1/4$ and the response radius specified in that
proposition. Both procedures recover $S$ with probability at least
$1-d^{-c'}-n^{-c''}$ if
\begin{equation}
\theta_{\min}\geq C\left\{
\sqrt{\frac{\log d}{n}}+
\frac{\sqrt{\log n}\log(\epsilon/\delta)}{n\epsilon}
\right\}.
\label{eq: sufficient beta-min}
\end{equation}
If additionally $|\varepsilon|\leq B_\varepsilon$ almost surely for a
public constant $B_\varepsilon$, replace $B_Y$ by
$B_xa_0+B_\varepsilon$ in the same penalty and sample-size conditions.
The finite-sample guarantee then holds for every $\alpha\in(0,1/4)$
under $\theta_{\min}\geq C_1\lambda$. In particular, taking
$\alpha=d^{-c'}<1/4$ gives recovery probability at least $1-d^{-c'}$ if
\begin{equation}
\theta_{\min}\geq C\left\{
\sqrt{\frac{\log d}{n}}+
\frac{\log(\epsilon/\delta)}{n\epsilon}
\right\}.
\label{eq: bounded noise}
\end{equation}
Here $C>0$ is sufficiently large and may depend only on the fixed public
model and tuning bounds and the displayed confidence exponents.
\end{theorem}

{\bf Buffer selection via certificate aggregation.}
The calibration of three buffer parameters
$(\psi^{\rm Param},\psi^{\rm Inv},\psi^{\rm Cross})$ involves a
natural trade-off: overly tight choices reduce the empirical buffer slack,
whereas overly conservative choices enlarge $\tau_{\rm primal}$ and
$\tau_{\rm dual}$ and weaken the normalized KKT margins. We therefore provide a practical selection rule allowing a public family of candidate buffer triples and aggregate their stability certificates rather than committing to a single calibration. By Proposition~\ref{prop:maximum-of-sub-distances}, this optimization does not require running multiple private procedures or composing privacy losses.  Specifically, let $\mathcal G\subset(0,\infty)^3$ be a nonempty family fixed independently of the protected dataset. For each $g=(\psi_g^{\rm Param},\psi_g^{\rm Inv},\psi_g^{\rm Cross})\in\mathcal G$ and $t\in\mathcal T=\{1,\ldots,n\}$, let $\gamma_{t,g}$ denote the score in \eqref{eq: adaptive candidate score} computed with sparsity index $t$ and buffer triple $g$. Define
\[
\gamma_{t,\mathcal G}(\bZ)
:=
\sup_{g\in\mathcal G}\gamma_{t,g}(\bZ),
\quad \mbox{ and } \quad
\gamma_{{\rm ad},\mathcal G}(\bZ)
:=
\max_{t\in\mathcal T}\gamma_{t,\mathcal G}(\bZ).
\]
\begin{corollary}
\label{cor:public-buffer-aggregation}
For every $t\in\mathcal T$, $\gamma_{t,\mathcal G}$ is a sub-distance for the common LASSO selector $\widehat S_\lambda$, and so is $\gamma_{{\rm ad},\mathcal G}$. Let $A=\widehat S_\lambda(\bZ)$ and $a=|A|$. If $A=\emptyset$, $a>n$, or $\widehat\Sigma_{AA}(\bZ)$ is singular, then $\gamma_{{\rm ad},\mathcal G}(\bZ)=0$; otherwise,
\[
\gamma_{{\rm ad},\mathcal G}(\bZ)
=
\sup_{g\in\mathcal G}\gamma_{a,g}(\bZ).
\]
Consequently, Saturated PTR can be applied once to the aggregated score without privacy composition over $\mathcal G$, and the exact oracle-equivalence identity in \eqref{eq:adaptive-pointwise-oracle-equivalence} continues to hold with $\gamma_s$ and $\gamma_{\rm ad}$ replaced by $\gamma_{s,\mathcal G}$ and $\gamma_{{\rm ad},\mathcal G}$, respectively.
\end{corollary}
In practice, one may take $\mathcal G$ to be a finite, publicly specified grid over conservative ranges, such that the algorithm is implementable and measurable. After a single LASSO fit, only the scalar stability scores need to be evaluated across the grid. By construction, enlarging the candidate set cannot decrease the aggregated score. Therefore, if $\mathcal G$ contains the calibration choice for the three buffer parameters $(\psi^{\rm Param},\psi^{\rm Inv},\psi^{\rm Cross})$ specified in Proposition~\ref{thm: sufficient beta-min}, both the fixed-sparsity and adaptive aggregated procedures inherit their respective fixed-buffer recovery guarantees stated in Proposition~\ref{thm: sufficient beta-min} and Theorem~\ref{cor: bounded noise rate}.

The penalty \(\lambda\) is likewise a public input, fixed from public information or a non-sensitive pilot sample. Choosing $\lambda$ from protected data requires a privacy guarantee for that tuning step and appropriate accounting for the resulting procedure.

\begin{remark}[Comparison with existing private support recovery methods]\label{rem:existing-private-methods}
The sufficient beta-min condition for the $(\epsilon,\delta)$-DP perturbation-stability LASSO of \citet{Guha-COLT-2013} contains a privacy term of order \(s^{3/2}/(n\epsilon)\), up to logarithmic factors, and its statistical recovery guarantee is stated at constant success probability. The pure-DP best subset selection (BSS) \citep{Roy-2024} and top-$R$ \citep{Prastakos-NIPS-2025} algorithms provide high-probability exact support recovery under their respective design regularity conditions, with private beta-min terms of orders $\sqrt{s\log(d)/(n\epsilon)}$ and $s\sqrt{\log(d)/(n\epsilon)}$, respectively. Both subset-selection approaches use known sparsity; exact BSS uses the exponential mechanism over all size-$s$ supports, while top-$R$ uses mixed-integer programming. The analyses of two other methods do not provide the worst-case privacy guarantee: the privacy guarantee of the mistakes procedure of \citet{Prastakos-NIPS-2025} requires a data-dependent gap condition, while the computationally efficient MCMC of \citet{Roy-2024} relies on a model-dependent rapid-mixing bound.

Under Assumption~\ref{assumption: regularity}, Theorem~\ref{cor: bounded noise rate} gives high-probability exact recovery for the computationally efficient Saturated LASSO with no explicit dependence on $s$ or $d$ in its private beta-min term. Its adaptive version has the same recovery probability without knowing $s$, with otherwise identical public tuning parameters. Complementary comparisons with Noisy IHT are provided in Section~\ref{sec:estimation-versus-recovery}.
\end{remark}

\section{The cost of privacy for exact support recovery}\label{sec: cost of privacy}

In this section, we investigate the information-theoretic limits of exact support recovery under $(\epsilon,\delta)$-DP. Section \ref{sec:minimax-lower-bound} establishes the minimax lower bound, characterizes its \(\delta\)-dependent phase transition, and examines attainability; Section \ref{sec:estimation-versus-recovery} contrasts exact support recovery with parameter estimation in terms of privacy cost.

\subsection{Minimax lower bound and \texorpdfstring{$\delta$}{delta}-dependent phase transition}\label{sec:minimax-lower-bound}

Let $\mathcal L(Z)$ denote the distribution of $Z=(Y,X^\top)^\top$. For fixed constants
$\sigma_x\geq1$, $\sigma>0$, and $0<c_\Sigma\leq1\leq C_\Sigma<\infty$, define the
$s$-sparse linear regression model class
\[
\begin{aligned}
\mathcal P(s,\theta_{\min}):=\bigg\{  \mathcal L(Z):\;
Y&=X^\top\btheta+\varepsilon,\
\btheta\in\Theta(s,\theta_{\min}),\ \varepsilon\perp X,  \ X\sim\mathrm{subGaussian}(0_d,\sigma_x^2),\\
&\
c_\Sigma I_d\preceq\cov(X)\preceq C_\Sigma I_d, \  \varepsilon\sim\mathrm{subGaussian}(0,\sigma^2)\bigg\}.
\end{aligned}
\]
The minimax risk for exact support recovery is
\[
    \mathcal R
    :=
    \inf_{\mathcal M\in\mathsf{DP}(\epsilon,\delta)}
    \sup_{P\in\mathcal P(s,\theta_{\min})}
    \mathbb P_{P,\mathcal M}
    \left(
    \mathcal M(\bZ)\neq\supp\{\btheta(P)\}
    \right),
\]
where $\bZ\sim P^{\otimes n}$ and the probability includes the randomness of the mechanism $\mathcal M$.

  The following theorem is obtained by constructing a grid of support-swap alternatives that identifies the hardest problem at each privacy level.
\begin{theorem} \label{thm: lower-bound}
Assume $d\ge 3s$, $s\ge 4$, $d-s+1\ge 16$, $0<\epsilon\le 1$, and $0\le\delta\le\epsilon/16$. If
\begin{equation}
 \theta_{\min}
 \le
 \sigma\left\{
 \frac{1}{4\sqrt2}
 \sqrt{\frac{\log(d-s+1)}{n}}
 +
 \frac{\sqrt2}{48n\epsilon}
 \max_{1\le r\le s}
 \frac{
 \min\{\log\left(\binom sr\binom{d-s}{r}\right),\log(\frac{\epsilon}{\delta})\}}{\sqrt r}
 \right\},
 \label{eq:complete-lower-bound}
\end{equation}
then $\mathcal R\ge\frac{1}{2}$.
\end{theorem}

\noindent{\bf Existing $\delta$-dependent lower bounds.} Lower bounds with explicit dependence on $\delta$ under $(\epsilon,\delta)$-DP have appeared in several settings \citep{Barber-Duchi-2014,Steinke-JPC-2016,cai2026cost}. These privacy bounds typically take the minimum of two terms: a global statistical complexity term and a $\delta$-dependent privacy  cutoff. In contrast, \eqref{eq:complete-lower-bound} applies the privacy cutoff separately to the entire family of $r$-swap alternatives and then optimizes over $r$. Consequently, varying $\delta$ changes not only the magnitude of the privacy cost but also the statistical scale of the least-favorable alternatives.

\begin{theorem}
\label{thm:phase-transition}
Assume the conditions on $d,s,\epsilon,\delta$ imposed in Theorem \ref{thm: lower-bound} hold. Also assume $ \log(e(d-s)/s)\ge5\log(es)$.  There exists a universal constant $c>0$ such that $\mathcal R\ge1/2$ if
\begin{equation}
 \theta_{\min}
 \le
 c\sigma\left\{
 \sqrt{\frac{\log(d-s+1)}{n}}
 +
 \theta_{\mathrm{priv}}
 \right\}, \quad \mbox{ where}
 \label{eq:three-phase-beta-min}
\end{equation}
\begin{numcases}{\theta_{\mathrm{priv}}=}
 \displaystyle
 \frac{1}{n\epsilon}\log(\frac{\epsilon}{\delta}),
 &
 if \ $\displaystyle
 \frac{\epsilon}{s(d-s)}
 \le\delta\le\frac{\epsilon}{16}$
 \label{eq:moderate-delta}
 \\[1mm]
 \displaystyle
 \frac{1}{n\epsilon}
 \sqrt{
 \log(\frac{\epsilon}{\delta})
 \log(\frac{e(d-s)}{s})},
 &
 if \ $\displaystyle
 \frac{\epsilon}{\binom{d-s}{s}}
 <\delta<
 \frac{\epsilon}{s(d-s)}$
 \label{eq:transition-phase}
 \\[1mm]
 \displaystyle
 \frac{\sqrt s}{n\epsilon}
 \log(\frac{e(d-s)}{s}),
 &
 if \ $\displaystyle
 0\le\delta\le
 \frac{\epsilon}{\binom{d-s}{s}}$
 \label{eq:pure-dp-regime}.
\end{numcases}

\end{theorem}
The additional condition $\log(e(d-s)/s)\ge5\log(es)$,
which holds for typical high-dimensional settings when \(d\) grows sufficiently rapidly relative to \(s\), gives the simple
expression in \eqref{eq:transition-phase}. The moderate-\(\delta\) rate in \eqref{eq:moderate-delta} follows directly from \eqref{eq:complete-lower-bound} and does not require this additional dimension condition.

When $\epsilon/(s(d-s))
\leq\delta\leq \epsilon/16$,
the maximum in (\ref{eq:complete-lower-bound}) is attained at $r=1$. The hardest comparison is therefore local: the true support is compared with alternatives obtained by a single swap. The private beta-min term is of order $\log(\epsilon/\delta)/(n\epsilon)$, and this interval defines the moderate-$\delta$ phase when $d$ and $s$ do not enter the private rate directly. As $\delta$ decreases below $\epsilon/\{s(d-s)\}$, the one-swap contribution saturates, while larger support-swap scales become progressively relevant, giving the transition rate in \eqref{eq:transition-phase}. Once $0\leq\delta\leq\epsilon /\binom{d-s}{s}$,
the full-swap choice $r=s$ yields the effectively pure-DP rate $\sqrt{s}\log(e(d-s)/s)/(n\epsilon)$. That is, $\delta$ is now too small to improve the lower-bound order over pure DP. This also explains why $d$ and $s$ reappear: the difficulty is driven by the global collection of possible supports rather than by a local swap.

We now compare the lower bound with the recovery guarantees  on a common model class.

\begin{corollary}[Minimax optimality of Saturated LASSO]
\label{cor:minimax-optimality}
Let $\mathcal P_{\rm reg}(s,\theta_{\min})$ be the subclass of
$\mathcal P(s,\theta_{\min})$ satisfying
Assumption~\ref{assumption: regularity} with fixed public regularity
constants, where $B_x>1$ and $a_1,a_2\geq1$, and let
$\mathcal P_{{\rm reg},{\rm bd}}(s,\theta_{\min})$ be its subclass with
$|\varepsilon|\leq\sigma$ almost surely. Under the dimensional and
privacy conditions of Theorem~\ref{thm: lower-bound}, both classes are
nonempty exactly when $s\theta_{\min}<a_0$. On these nonempty classes,
the lower-bound rate in \eqref{eq:complete-lower-bound} also holds up to universal constants.
The corresponding recovery guarantees
\eqref{eq: sufficient beta-min} and \eqref{eq: bounded noise} hold
uniformly under the conditions of Theorem~\ref{cor: bounded noise rate},
with $B_\varepsilon=\sigma$ in the bounded-noise case.

Consequently, under these recovery conditions, whenever the
sufficient signal threshold is smaller than $a_0/s$, both fixed-sparsity
and adaptive Saturated LASSO are minimax optimal on these classes up to
logarithmic factors in $n$ and $1/\delta$, uniformly over
$0<\delta\leq\epsilon/16$.
When
$\epsilon/\{s(d-s)\}\leq\delta\leq\epsilon/16$, the rates match up to
the $\sqrt{\log n}$ response-clipping factor under sub-Gaussian noise
and up to constants under bounded noise.

\end{corollary}

 The sharper moderate-$\delta$ conclusion applies to polynomially small privacy parameters in a broad range of high-dimensional recovery settings. For example, take $\delta=n^{-(1+c)}$ for fixed $c>0$ \citep{Cai-Wang-Zhang-2021} and fixed $0<\epsilon\leq1$. The moderate-$\delta$ condition requires $s(d-s)\geq\epsilon n^{1+c}$. Choosing $d=\lceil n^{2+c}\rceil$ and $s=\lfloor n^{1/4}\rfloor$ satisfies this condition, the sample-size requirements of Theorem~\ref{cor: bounded noise rate}, and a nonempty recovery window for all sufficiently large $n$, with the public model and tuning bounds held fixed.

\paragraph{Attainability under independent Gaussian design.}
A complementary achievability result, including detailed constructions and proofs, is presented in Section~\textcolor{Blue}{S6} of the
supplementary material.
Under the conditions on $d,s,\epsilon,\delta$ in Theorem~\ref{thm:phase-transition},
Algorithm~\textcolor{Blue}{S3} satisfies worst-case $(\epsilon,\delta)$-DP.
Consider model~\eqref{eq: linear regression model} with $X\sim N(0,\sigma_x^2 I_d)$,
independent noise $\varepsilon\sim N(0,\sigma^2)$, known
sparsity $s$, and $\|\btheta^*\|_2\le R$.
Corollary~\textcolor{Blue}{S2} shows that, for every fixed $c>0$, there exists
a constant $C_c>0$ depending only on $c$ such that,  if
\[
\theta_{\min}\ge
C_c\sqrt{R^2+\frac{\sigma^2}{\sigma_x^2}}
\left\{
\sqrt{\frac{\log d}{n}}+\theta_{\mathrm{priv}}
\right\},
\]
then $\widehat S_{\mathrm{priv}}$ recovers the true support $S$ with probability at least $1-d^{-c}$,
where $\widehat S_{\mathrm{priv}}$ is the output of
Algorithm~\textcolor{Blue}{S3} and $\theta_{\mathrm{priv}}$ is given
in \eqref{eq:moderate-delta}--\eqref{eq:pure-dp-regime}.
This Gaussian model is nonempty exactly when $\sqrt{s}\theta_{\min}\leq R$, so the displayed sufficient threshold must not exceed $R/\sqrt{s}$ to give a nonempty recovery range.
Therefore, for fixed positive $R,\sigma,$ and $\sigma_x$, it matches the three lower-bound orders in $n,d,s,\epsilon,$ and $\delta$ within this Gaussian model. Algorithm~\textcolor{Blue}{S3} has exponential computational complexity and is an information-theoretic complement to the computationally efficient, sparsity-adaptive Saturated LASSO.

\subsection{Support recovery versus coefficient estimation}
\label{sec:estimation-versus-recovery}

We compare established recovery and existing estimation guarantees on related but nonidentical distribution classes, with model-scale and regularity constants held fixed. This rate-level comparison highlights a potential difference in how privacy affects the two problems,  and provides a useful geometric interpretation of this distinction.

Consider a parameter with $s$ nonzero coordinates of a common magnitude
$\theta_{\min}$. In the non-private setting, the exact-recovery scale
$\theta_{\min}\asymp \sigma\sqrt{\log d/n}$ corresponds to the
$\ell_2$ signal scale
$\sqrt{s}\theta_{\min}\asymp \sigma\sqrt{s\log d/n}$.
When $\log(ed/s)\asymp\log d$, this has the same order as the global
minimax $\ell_2$ estimation rate over $s$-sparse parameter classes;
see \citet{Ndaoud-ICALT-2019}.

Under $(\epsilon,\delta)$-DP, Corollary~\ref{cor:minimax-optimality} shows that Saturated LASSO attains the moderate-$\delta$ recovery rate under bounded noise. The corresponding $\ell_2$ signal threshold has a privacy term of order $\sigma\sqrt{s}\log(\epsilon/\delta)/(n\epsilon)$; under sub-Gaussian noise, the recovery guarantee incurs the additional $\sqrt{\log n}$ response-clipping factor. For comparison, the available $\ell_2$-error upper bound for noisy iterative hard thresholding (Noisy IHT) has an $s\log d/(n\epsilon)$ privacy factor, together with clipping and logarithmic factors \citep{Cai-Wang-Zhang-2021,cai2026cost}. The recovery signal threshold thus has a more favorable explicit dependence on $s$ than the corresponding term in this estimation bound.

The geometry of the targets motivates different privacy protection strategies.
Noisy IHT perturbs numerical quantities during optimization and uses
noisy top-$k$ selection, which allocates privacy across
coordinate comparisons.
Saturated LASSO instead uses its internal fit to certify a discrete
support decision and applies the private release rule once.
When the certificate reaches the saturation threshold established in Proposition \ref{prop:finite-saturation}, Saturated PTR
releases the selected support with probability one, without releasing
coefficient estimates. Section~\ref{sec: simulation} provides a numerical comparison between Saturated LASSO and Noisy IHT in terms of exact support recovery performance.
Section~\textcolor{Blue}{S1} of the supplementary material illustrates this
difference through sufficient recovery conditions under an independent
design.  As \(\delta\) decreases, larger support swaps become relevant and the private recovery boundary acquires dependence on \(s\) and \(d\), as shown in \eqref{eq:transition-phase}–\eqref{eq:pure-dp-regime}.

\section{Numerical results}\label{sec: numerical results}

This section evaluates the numerical performance of Saturated LASSO (Algorithm~\ref{alg:saturated-lasso}) and illustrates  Adaptive Saturated LASSO (Algorithm~\ref{alg:adaptive-saturated-lasso}).
Section~\ref{sec: simulation} examines the cost of privacy and compares the method with existing private procedures.
Theorem~\ref{thm:adaptive-oracle-equivalence} establishes identical finite-sample exact-recovery risks for Adaptive Saturated LASSO and its oracle fixed-$s$ counterpart under common public tuning.
Section~\ref{sec: real-data} illustrates the adaptive procedure using public American Community Survey housing covariates, and provides a common setting for examining support selection, certification, and sensitivity to public parameters.
Throughout this section, selection accuracy is the proportion of correct support recoveries over $B=500$ replications.
For a private procedure, this is the frequency of correct release; both incorrect release and abstention count as recovery failures.
Error bars are pointwise $95\%$ Wilson score intervals \citep{Brown-Cai-DasGupta-2001}.
For each fixed configuration, observing zero incorrect releases in \(500\) independent replications yields a one-sided \(95\%\) exact upper confidence bound of approximately \(0.60\%\) for the unconditional probability of an incorrect release.

\subsection{Simulation studies}\label{sec: simulation}
The data $(X_1,Y_1),\ldots,(X_n,Y_n)$  are i.i.d. from \eqref{eq: linear regression model}. Given  $s$ and  $\theta_{\min}$, we sample $S$ uniformly without replacement, set $\theta_j^*=\theta_{\min}$ for $j\in S$, and set the remaining coefficients to zero; this parameter is fixed across replications within a setting. We generate $V\sim N(0,\Sigma_V)$, where $(\Sigma_V)_{ij}=\sigma_0\rho^{|i-j|}$, and set $X_j=Q_{\sigma_x}(V_j)$ for $j\in[d]$. The noise is  $N(0,\sigma^2)$. We set $B_Y = \sqrt{\log(n)}$ for Saturated LASSO, and  $\delta = 10n^{-1.1}$ as specified in
\cite{Chang-JASA-2026}.
Tuning parameters were calibrated using independent synthetic pilot samples and frozen before evaluation.

\begin{figure}[!ht]
    \centering
    \includegraphics[width=1.0\linewidth]{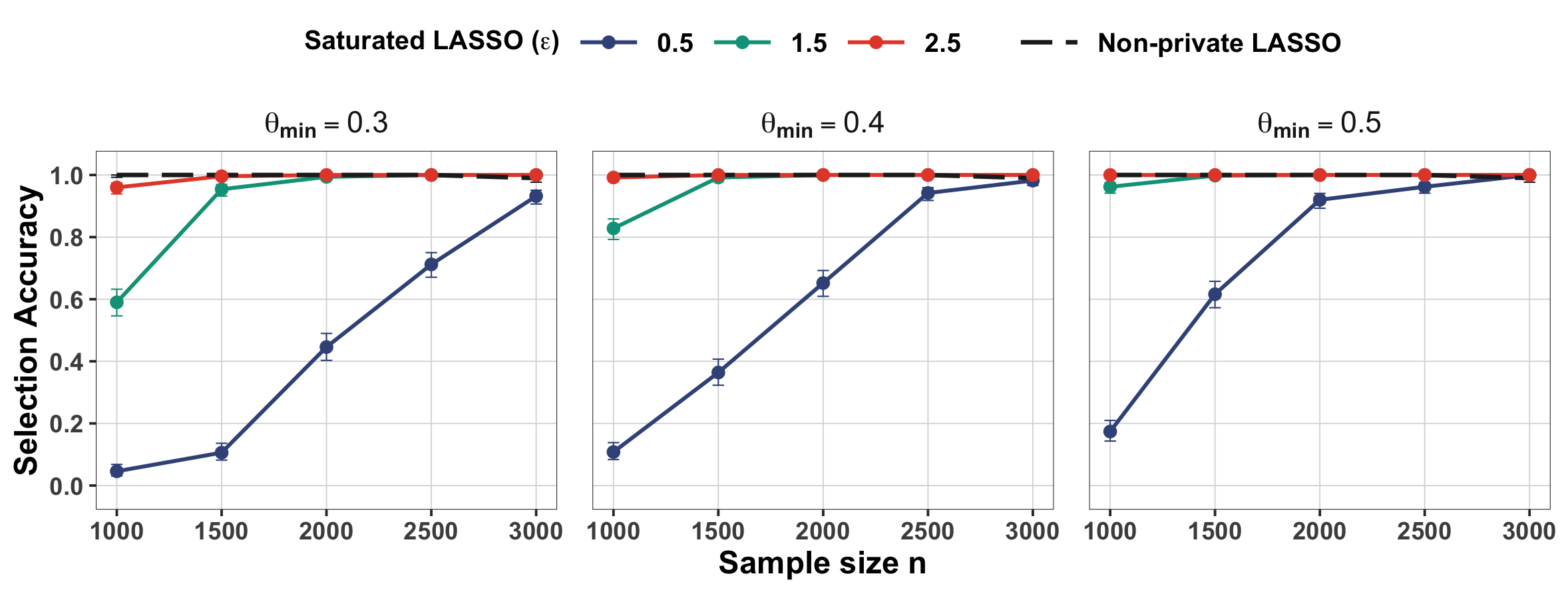}
    \caption{Estimated selection accuracy for support recovery of Saturated LASSO over $B = 500$ replications across different sample sizes $n\in \{1000,1500,2000,2500,3000\}$, privacy parameters $\epsilon \in \{0.5,1.5,2.5\}$, and minimum signal strengths $\theta_{\min}\in\{0.3,0.4,0.5\}$.  We set $\sigma_x = \sigma_0 = 1$, $\sigma = 0.5$, $\rho = 0.2$, $s = 5$ and $d = 2n$. }
    \label{fig:cost of privacy}
\end{figure}

\Needspace{3\baselineskip}
\noindent{\bf The cost of privacy for support recovery.}  Figure~\ref{fig:cost of privacy} illustrates the cost of privacy for exact support recovery by evaluating Saturated LASSO (Algorithm~\ref{alg:saturated-lasso}) across sample sizes \(n\in\{1000,1500,2000,2500,3000\}\), privacy budgets \(\epsilon\in\{0.5,1.5,2.5\}\), and minimum signal strengths \(\theta_{\min}\in\{0.3,0.4,0.5\}\).
We set \(s=5\) and \(d=2n\), and include the non-private LASSO as a benchmark.
The private procedure and non-private benchmark use separately calibrated penalties, with calibration specific to the simulated configuration.
Under the strongest privacy constraint, \(\epsilon=0.5\), recovery improved substantially with both sample size and signal strength: for example, when \(\theta_{\min}=0.3\), accuracy increased from \(0.046\) at \(n=1000\) to \(0.932\) at \(n=3000\).
The candidate support entering Saturated LASSO was correct in every evaluation run, so its observed recovery failures arose from abstention.
The privacy cost became much smaller as \(\epsilon\) increased.
With \(\epsilon=1.5\), nearly perfect recovery was attained for all signal strengths by \(n=2000\), while \(\epsilon=2.5\) performed nearly identically to the non-private benchmark.

\Needspace{3\baselineskip}
\noindent{\bf Performance comparison of the algorithms.}
 We next examine differences between private procedures as sparsity increases, across different sample sizes and signal strengths.
We compare Saturated LASSO with perturbation-stability LASSO \citep{Guha-COLT-2013} (Stability LASSO) and Noisy IHT \citep{Cai-Wang-Zhang-2021}.
The former is the closest benchmark for stability-based private release of a LASSO support, whereas the latter examines support recovery obtained from private coefficient estimation, corresponding to the distinction developed in Section~\ref{sec:estimation-versus-recovery}.
The original perturbation-stability LASSO uses add/remove adjacency; we adapt it to record-replacement adjacency so that all implemented procedures satisfy the same prescribed \((\epsilon,\delta)\)-DP requirement.
For Noisy IHT, we select the indices of the \(s\) largest entries, in absolute value, of its final coefficient estimate.
The experiments focus on these two methodological comparisons.
Complementary subset-selection approaches, including the pure-DP BSS and top-\(R\) procedures of \cite{Roy-2024} and \cite{Prastakos-NIPS-2025}, are discussed in Remark~\ref{rem:existing-private-methods}.
Implementation and calibration details for both simulation studies, including the record-replacement privacy argument for perturbation-stability LASSO, are provided in Section~\textcolor{Blue}{S7} of the supplement.

\begin{figure}[!ht]
    \centering
    \includegraphics[width=1.0\linewidth]{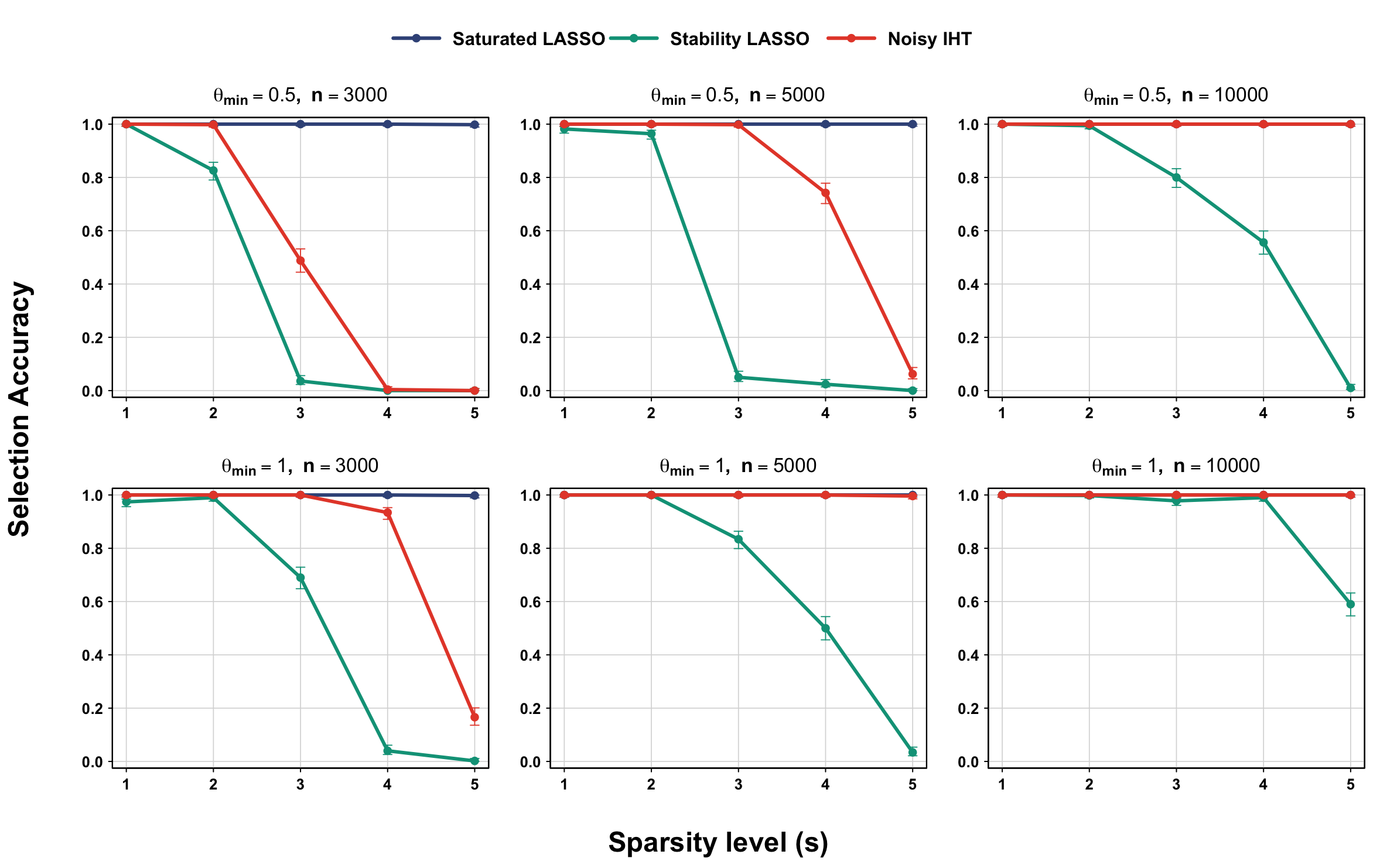}
    \caption{Estimated selection accuracy for support recovery of the three algorithms over $B = 500$ replications with different sample sizes $n\in\{3000,5000,10000\}$, sparsity levels $s\in\{1,\ldots,5\}$, and minimum signal strengths $\theta_{\min}\in\{0.5,1.0\}$.  We set $\sigma_x = \sigma_0 = 1$, $\sigma = 0.5$, $\rho = 0.1$, and $\epsilon = 4$. }
    \label{fig:recovery comparison}
\end{figure}

Figure~\ref{fig:recovery comparison} reports the estimated selection accuracy of the three algorithms across sample sizes $n\in\{3000,5000,10000\}$, sparsity levels $s\in\{1,\ldots,5\}$,  minimum signal strengths $\theta_{\min}\in\{0.5,1.0\}$, and $d = n$.
Saturated LASSO achieved nearly perfect recovery in each of the $30$ configurations. All three procedures achieved this accuracy at $s=1$, but their performance separated as sparsity increased. For example, at $n=3000$, $\theta_{\min}=0.5$, and $s=3$, the recovery rates were $1.000$, $0.036$, and $0.488$ for Saturated LASSO, perturbation-stability LASSO, and Noisy IHT, respectively. Larger samples and stronger signals generally improved the performance of the competing procedures; in particular, Noisy IHT matched Saturated LASSO's observed recovery across all sparsity levels at $n=10000$ for both signal strengths.

To interpret the gap between the two LASSO-based procedures, we distinguish support selection from private release.
The simulation outputs show that both procedures selected the true support before release in all $15{,}000$ evaluation runs per procedure.
Thus, the lower recovery rates of perturbation-stability LASSO arose entirely from withholding correct selections.
At $n=3000$, $\theta_{\min}=0.5$, and $s=3$, it withheld release in $96.4\%$ of replications, whereas Noisy IHT always returned a support and returned an incorrect support in $51.2\%$.
These results are consistent with the combined effect of the KKT certificate and the saturated release rule, as well as the more favorable sparsity dependence discussed in Remark~\ref{rem:existing-private-methods}.
The comparison with Noisy IHT provides an illustration of the distinction developed in Section~\ref{sec:estimation-versus-recovery}: guarantees for private coefficient estimation do not by themselves ensure comparable exact support recovery.

\subsection{Semi-synthetic case study}\label{sec: real-data}
We use covariates from the 2023 American Community Survey (ACS), one-year California housing sample, to illustrate Adaptive Saturated LASSO under an observed covariate distribution. We examine how support selection and stability certification jointly shape private recovery as study size and privacy parameters vary. The data contain \(N=140{,}143\) occupied housing units and yield \(d=54\) candidate terms.

We retain the observed covariates, using fixed public transformations to ensure $\|X_i\|_\infty\leq1$, and generate responses according to model~\eqref{eq: linear regression model}.
The active terms are the number of rooms, household size, and indicators of electricity as the main heating fuel and construction before 1940.
Each transformed active term has coefficient $0.125$ and all other coefficients are zero, giving $s=4$ and $\theta_{\min}=0.125$.
Independent noise uniform on $[-0.5,0.5]$ gives $|Y_i|\leq1$.
This construction preserves the empirical covariate dependence while providing known ground truth.
Because the ACS microdata are public and the responses are simulated, the experiment diagnoses a hypothetical private analysis rather than analyzing confidential outcomes.

We apply Adaptive Saturated LASSO  with the public buffer aggregation of Corollary~\ref{cor:public-buffer-aggregation} and use the corresponding non-private LASSO as a benchmark.
We fix $\lambda=0.008$, $B_X=B_Y=1$, and $\delta=5\times10^{-6}$, and examine sample sizes from $500$ to $140{,}143$ and privacy budgets $\epsilon\in\{0.5,1,2,4\}$.
These tuning values and the buffer grid are fixed across replications.
At each study size, the same LASSO fit and certificate are used across privacy budgets, so varying $\epsilon$ changes the exact recurrence probability $q_{\lceil\gamma\rceil}$ without changing the candidate support.
Each setting uses $B=500$ independent replications.
Section~\textcolor{Blue}{S8} of the supplement provides data and implementation details, together with the supporting diagnostics and sensitivity analyses of the privacy parameters $(\epsilon,\delta)$ and penalty level $\lambda$.

\begin{figure}[!ht]
  \centering
  \setstretch{1}\setlength{\parskip}{0pt}
  \includegraphics[width=0.70\linewidth]{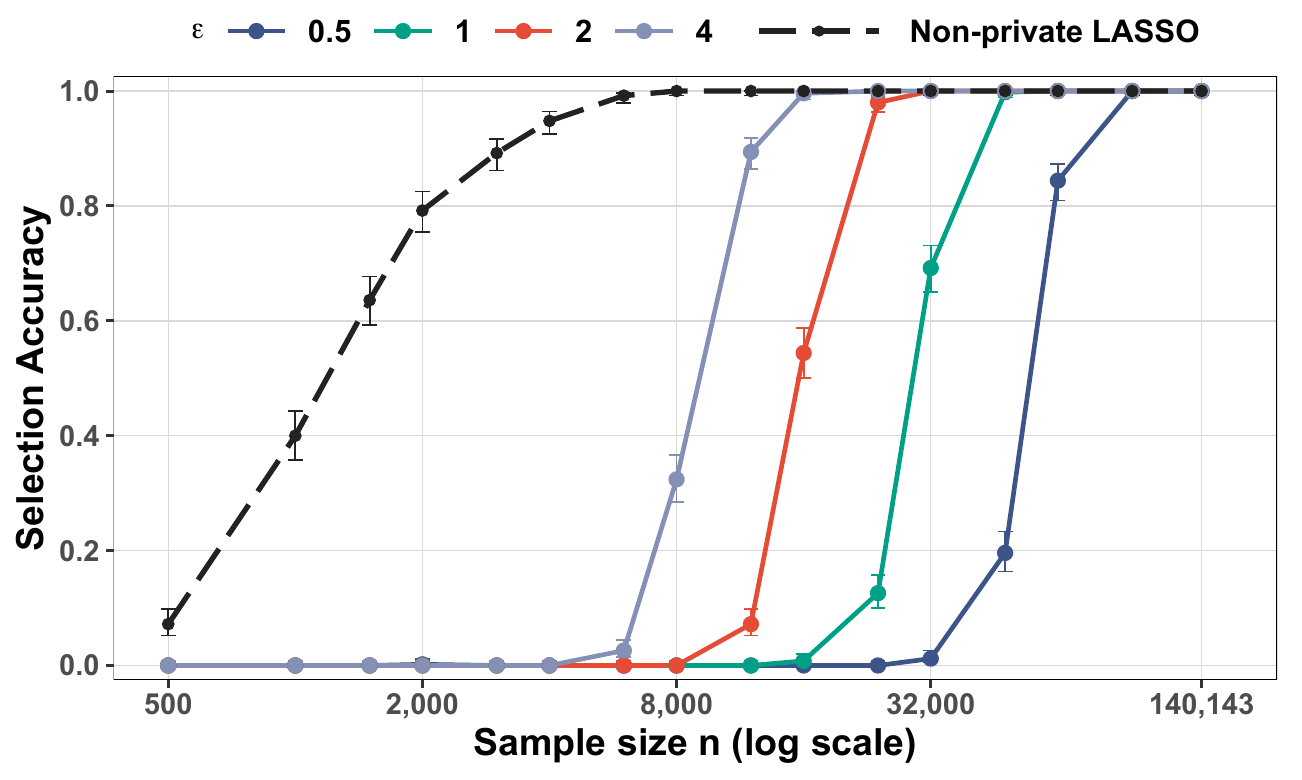}
  \caption{Estimated selection accuracy over $B=500$ replications, with $d=54$, $s=4$, $\theta_{\min}=0.125$, $\lambda=0.008$, and $\delta=5\times10^{-6}$.
Colored solid curves with dots give Adaptive Saturated LASSO accuracies at $\epsilon\in\{0.5,1,2,4\}$; the black dashed curve with dots gives the non-private LASSO accuracy.}
  \label{fig:acs-recovery}
\end{figure}
Figure~\ref{fig:acs-recovery} shows that non-private recovery reached high accuracy at smaller study sizes than private recovery.
Non-private accuracy increased from $40.0\%$ at $n=1000$ to $100\%$ at every evaluated size $n\geq8000$.
At $\epsilon=2$, private accuracy increased from $54.4\%$ at $n=16{,}000$ to $98.0\%$ at $n=24{,}000$.
Increasing the privacy budget also narrowed the gap: at $n=32{,}000$, private accuracy increased from $69.2\%$ at $\epsilon=1$ to $100\%$ at $\epsilon=2$. No incorrect release was observed at the settings displayed in Figure~\ref{fig:acs-recovery}, so all observed private recovery failures there were abstentions.
A complementary sensitivity analysis of the penalty level at $n=32{,}000$ shows that private recovery
remained robust across the tested penalties at mild privacy budgets $\epsilon\in\{2,3,4\}$, and varied markedly at the more stringent privacy level $\epsilon=1$, with non-private recovery unchanged (Supplementary Figure~\textcolor{Blue}{S3}).

The case study distinguishes the accuracy of the selected support from the stability certification governing its private release. The certificate measures robustness to record replacement and does not, by itself, establish that the selected support is correct. At smaller study sizes, selected variables vary across samples, whereas withholding can persist even after non-private recovery becomes accurate. The KKT construction in Section~\ref{sec: perturbation and stability} ties certification to separation from support-changing boundaries relative to record-replacement perturbations. Consequently, weak primal or dual separation may also indicate fragility of the underlying LASSO selection, independently of the privacy mechanism. Additional observations can both improve selection accuracy and strengthen the stability certificate, thereby increasing the frequency of private release at a fixed privacy level. Increasing \(\epsilon\) can likewise increase the release frequency, but it leaves the fitted model and its certificate unchanged and therefore provides no additional statistical evidence for the selected variables. Following abstention, investigators may instead collect additional data before reconsidering private release.

\section{Discussion}\label{sec: discussion}

This paper highlights how the discreteness of a statistical target can shape the mechanism design and cost of privacy protection. A selected support can remain unchanged under data perturbations that may otherwise change the estimated parameters, allowing support selection and private release to be treated separately. In Saturated LASSO, statistical accuracy depends on signal strength and  regularity conditions, whereas privacy holds on every dataset; once the stability certificate reaches the finite saturation threshold $K_{\epsilon,\delta}$, a correct selection is released with probability one. Sparsity adaptation exploits the same separation by aggregating certificates for a common support decision without repeated private releases. The minimax analysis provides a complementary perspective: as \(\delta\) decreases, increasingly large support changes govern the privacy cost, connecting local alternatives in the moderate-\(\delta\) regime to the global collection of competing supports near effectively pure privacy regime. On the regular model classes considered here and under the stated recovery conditions, uniform logarithmic-factor optimality and sharper moderate-\(\delta\) matching establish the effectiveness of Saturated LASSO. Attaining all three lower-bound orders with an efficient procedure under dependent designs remains open.

An immediate direction for future research is to extend stability-certified release beyond the LASSO. Many statistical procedures produce discrete structures that remain unchanged over neighborhoods of the data, including grouped-variable selectors, generalized linear model selectors, change-point estimators, clustering procedures, and graph-selection methods. The central challenge is to construct stability certificates that are both computationally tractable and sufficiently sharp to preserve the statistical performance of the underlying non-private procedure. For example, in sparse precision-matrix estimation, such certificates could enable private release of conditional-dependence graphs arising in genomics and network analysis. More broadly, an important theoretical direction is to characterize general geometric conditions under which a discrete statistical procedure admits an efficiently computable sub-distance certificate.

Exact support recovery may also be too stringent under strong predictor dependence. When several correlated variables carry nearly interchangeable information, releasing groups of variables or a selected set with false discovery control may be more scientifically meaningful. Developing stability certificates for these weaker targets could produce useful procedures in settings where exact coordinatewise identification is impossible. Another related question is how to combine certification with post-selection estimation or inference under a shared privacy budget. The abstention outcome also merits further study as a formal decision option, including how its probability should be reported and incorporated into downstream scientific workflows.

\bibliographystyle{apalike}
\bibliography{ref}

@article{Wainwright-TIT-LASSO-2009,
  title   = {Sharp Thresholds for High-Dimensional and Noisy Sparsity Recovery Using {$\ell_1$}-Constrained Quadratic Programming ({L}asso)},
  author  = {Martin J. Wainwright},
  journal = {IEEE Transactions on Information Theory},
  year    = {2009},
  volume  = {55},
  number  = {5},
  pages   = {2183--2202},
  url     = {https://arxiv.org/abs/math/0605740}
}

@misc{Shen-arxiv-2026,
  title         = {Efficient Propose-Test-Release for Optimal Differentially Private Estimation},
  author        = {Tao Shen and Xin T. Tong and Wanjie Wang},
  year          = {2026},
  eprint        = {2605.03264},
  archivePrefix = {arXiv},
  primaryClass  = {stat.ME},
  note          = {arXiv:2605.03264},
  url           = {https://arxiv.org/abs/2605.03264}
}

@inproceedings{Dwork-Lei-STOC-2009,
  author    = {Cynthia Dwork and Jing Lei},
  title     = {Differential Privacy and Robust Statistics},
  booktitle = {Proceedings of the Forty-First Annual ACM Symposium on Theory of Computing},
  pages     = {371--380},
  year      = {2009},
}

@inproceedings{Awan-NeurIPS-2018,
  author    = {Jordan Awan and Aleksandra Slavkovi{\'c}},
  title     = {Differentially Private Uniformly Most Powerful Tests for Binomial Data},
  booktitle = {Advances in Neural Information Processing Systems},
  volume    = {31},
  pages     = {4208--4218},
  year      = {2018},
  url       = {https://papers.nips.cc/paper_files/paper/2018/hash/296472c9542ad4d4788d543508116cbc-Abstract.html}
}

@article{Desfontaines-POPETS-2022,
  author  = {Damien Desfontaines and James Voss and Bryant Gipson and Chinmoy Mandayam},
  title   = {Differentially Private Partition Selection},
  journal = {Proceedings on Privacy Enhancing Technologies},
  year    = {2022},
  volume  = {2022},
  number  = {1},
  pages   = {339--352},
}

@article{Cai-Wang-Zhang-2021,
  author  = {Cai, T. Tony and Wang, Yichen and Zhang, Linjun},
  title   = {The Cost of Privacy: Optimal Rates of Convergence for Parameter Estimation with Differential Privacy},
  journal = {The Annals of Statistics},
  year    = {2021},
  volume  = {49},
  number  = {5},
  pages   = {2825--2850},
}

@inproceedings{Roy-2024,
  title     = {On the Computational Complexity of Private High-Dimensional Model Selection},
  author    = {Roy, Saptarshi and Wang, Zehua and Tewari, Ambuj},
  booktitle = {Advances in Neural Information Processing Systems},
  volume    = {37},
  pages     = {36683--36716},
  year      = {2024},
  url       = {https://proceedings.neurips.cc/paper_files/paper/2024/hash/40889e2aedb47e64d42dd8b8de2fab4e-Abstract-Conference.html}
}

@misc{Cai-2023,
  title         = {Score Attack: A Lower Bound Technique for Optimal Differentially Private Learning},
  author        = {T. Tony Cai and Yichen Wang and Linjun Zhang},
  year          = {2023},
  eprint        = {2303.07152},
  archivePrefix = {arXiv},
  primaryClass  = {math.ST},
  note          = {arXiv:2303.07152},
  url           = {https://arxiv.org/abs/2303.07152}
}

@inproceedings{Dwork-2006,
  author    = {Dwork, Cynthia and McSherry, Frank and Nissim, Kobbi and Smith, Adam},
  title     = {Calibrating Noise to Sensitivity in Private Data Analysis},
  booktitle = {Theory of Cryptography: Third Theory of Cryptography Conference, {TCC} 2006},
  series    = {Lecture Notes in Computer Science},
  volume    = {3876},
  pages     = {265--284},
  publisher = {Springer},
  year      = {2006},
  isbn      = {978-3-540-32731-8}
}

@book{Dwork-2014,
  title={The Algorithmic Foundations of Differential Privacy},
  author={Dwork, Cynthia and Roth, Aaron},
  year={2014},
  publisher={Now Publishers Inc},
  series={Foundations and Trends in Theoretical Computer Science},
  volume={9},
  pages={211--407},
}

@article{Butucea-AoS-2018,
  title     = {Variable selection with {Hamming} loss},
  author    = {Butucea, Cristina and Ndaoud, Mohamed and Stepanova, Natalia A. and Tsybakov, Alexandre B.},
  journal   = {The Annals of Statistics},
  volume    = {46},
  number    = {5},
  pages     = {1837--1875},
  year      = {2018},
  publisher = {Institute of Mathematical Statistics},
}

@inproceedings{Steinke-FOCS-2017,
  author    = {Thomas Steinke and Jonathan Ullman},
  title     = {Tight Lower Bounds for Differentially Private Selection},
  booktitle = {2017 {IEEE} 58th Annual Symposium on Foundations of Computer Science (FOCS)},
  pages     = {552--563},
  year      = {2017},
  url       = {https://arxiv.org/abs/1704.03024},
}

@article{Steinke-JPC-2016,
  author  = {Thomas Steinke and Jonathan Ullman},
  title   = {Between Pure and Approximate Differential Privacy},
  journal = {Journal of Privacy and Confidentiality},
  year    = {2017},
  volume  = {7},
  number  = {2},
  pages   = {3--22},
  url     = {https://journalprivacyconfidentiality.org/index.php/jpc/article/view/648}
}

@inproceedings{Bafna-COLT-2017,
  author    = {Mitali Bafna and Jonathan Ullman},
  title     = {The Price of Selection in Differential Privacy},
  booktitle = {Proceedings of the 2017 Conference on Learning Theory},
  series    = {Proceedings of Machine Learning Research},
  volume    = {65},
  pages     = {151--168},
  year      = {2017},
  publisher = {PMLR},
  url       = {https://proceedings.mlr.press/v65/bafna17a.html}
}

@article{Wainwright-TIT-2009,
    author = {Martin J. Wainwright},
    title = {Information-Theoretic Limits on Sparsity Recovery in the High-Dimensional and Noisy Setting},
    journal = {IEEE Transactions on Information Theory},
    year = {2009},
    volume = {55},
    number = {12},
    pages = {5728--5741},
}

@article{Brown-Cai-DasGupta-2001,
  author  = {Brown, Lawrence D. and Cai, T. Tony and DasGupta, Anirban},
  title   = {Interval Estimation for a Binomial Proportion},
  journal = {Statistical Science},
  year    = {2001},
  volume  = {16},
  number  = {2},
  pages   = {101--133},
  doi     = {10.1214/ss/1009213286}
}

@inproceedings{Guha-COLT-2013,
  title     = {Differentially Private Feature Selection via Stability Arguments, and the Robustness of the Lasso},
  author    = {Thakurta, Abhradeep Guha and Smith, Adam},
  booktitle = {Proceedings of the 26th Annual Conference on Learning Theory},
  series    = {Proceedings of Machine Learning Research},
  volume    = {30},
  pages     = {819--850},
  year      = {2013},
  month     = {12--14 Jun},
  url       = {https://proceedings.mlr.press/v30/Guha13.html}
}

@article{Zhao-JMLR-2006,
  title   = {On Model Selection Consistency of Lasso},
  author  = {Peng Zhao and Bin Yu},
  journal = {Journal of Machine Learning Research},
  year    = {2006},
  volume  = {7},
  number  = {90},
  pages   = {2541--2563},
  url     = {https://jmlr.org/papers/v7/zhao06a.html}
}

@article{Zhang-JMLR-2009,
  title   = {On the Consistency of Feature Selection using Greedy Least Squares Regression},
  author  = {Tong Zhang},
  journal = {Journal of Machine Learning Research},
  volume  = {10},
  number  = {19},
  pages   = {555--568},
  year    = {2009}
}

@article{Fan-JASA-2001,
  author = {Fan, Jianqing and Li, Runze},
  title = {Variable Selection via Nonconcave Penalized Likelihood and its Oracle Properties},
  journal = {Journal of the American Statistical Association},
  year = {2001},
  volume = {96},
  number = {456},
  pages = {1348--1360},
}

@article{Zhang-AOS-2010,
  author = {Zhang, Cun-Hui},
  title = {Nearly Unbiased Variable Selection Under Minimax Concave Penalty},
  journal = {The Annals of Statistics},
  year = {2010},
  volume = {38},
  number = {2},
  pages = {894--942},
  url = {https://arxiv.org/abs/1002.4734},
}

@article{Cai-TIT-2011,
  title={Orthogonal Matching Pursuit for Sparse Signal Recovery With Noise},
  author={Cai, T. T. and Wang, Lie},
  journal={IEEE Transactions on Information Theory},
  year={2011},
  volume={57},
  number={7},
  pages={4680--4688},
}

@article{Cai-arxiv-2025,
  title   = {Differentially Private Estimation and Inference in High-Dimensional Regression with {FDR} Control},
  author  = {Zhanrui Cai and Sai Li and Xintao Xia and Linjun Zhang},
  journal = {Journal of Machine Learning Research},
  year    = {2026},
  volume  = {27},
  number  = {81},
  pages   = {1--54},
  url     = {https://www.jmlr.org/papers/v27/24-1413.html}
}

@inproceedings{Ndaoud-ICALT-2019,
  title     = {Interplay of Minimax Estimation and Minimax Support Recovery under Sparsity},
  author    = {Ndaoud, Mohamed},
  booktitle = {Proceedings of the 30th International Conference on Algorithmic Learning Theory},
  series    = {Proceedings of Machine Learning Research},
  volume    = {98},
  pages     = {647--668},
  year      = {2019},
  month     = {22--24 Mar},
  publisher = {PMLR},
  url       = {https://proceedings.mlr.press/v98/ndaoud19a.html}
}

@article{Cummings-HDSR-2024,
  author  = {Cummings, Rachel and Desfontaines, Damien and Evans, David and Geambasu, Roxana and Huang, Yangsibo and Jagielski, Matthew and Kairouz, Peter and Kamath, Gautam and Oh, Sewoong and Ohrimenko, Olga and Papernot, Nicolas and Rogers, Ryan and Shen, Milan and Song, Shuang and Su, Weijie and Terzis, Andreas and Thakurta, Abhradeep and Vassilvitskii, Sergei and Wang, Yu-Xiang and Xiong, Li and Yekhanin, Sergey and Yu, Da and Zhang, Huanyu and Zhang, Wanrong},
  title   = {Advancing Differential Privacy: Where We Are Now and Future Directions for Real-World Deployment},
  journal = {Harvard Data Science Review},
  year    = {2024},
  month   = {January},
  volume  = {6},
  number  = {1},
  url     = {https://hdsr.mitpress.mit.edu/pub/sl9we8gh/release/3}
}

@article{Lei-JRSSA-2018,
  author  = {Jing Lei and Anne-Sophie Charest and Aleksandra Slavkovi{\'c} and Adam Smith and Stephen Fienberg},
  title   = {Differentially Private Model Selection with Penalized and Constrained Likelihood},
  journal = {Journal of the Royal Statistical Society Series A: Statistics in Society},
  year    = {2018},
  volume  = {181},
  number  = {3},
  pages   = {609--633},
  url     = {https://academic.oup.com/jrsssa/article/181/3/609/7072010}
}

@INPROCEEDINGS{Reza-IEEESP-2017,
  author={Shokri, Reza and Stronati, Marco and Song, Congzheng and Shmatikov, Vitaly},
  booktitle={IEEE Symposium on Security and Privacy}, 
  title={Membership Inference Attacks Against Machine Learning Models}, 
  year={2017},
  volume={},
  number={},
  pages={3-18},
}

@inproceedings{Prastakos-NIPS-2025,
  title     = {Differentially Private High-Dimensional Variable Selection via Integer Programming},
  author    = {Petros Prastakos and Kayhan Behdin and Rahul Mazumder},
  booktitle = {Advances in Neural Information Processing Systems},
  volume    = {38},
  pages     = {129142--129182},
  year      = {2025},
  url       = {https://papers.nips.cc/paper_files/paper/2025/hash/a9401a4f93a42db4fccf935b7bad2362-Abstract-Conference.html}
}

@misc{Gu-DP-supportrecovery,
  title  = {Supplementary Material for ``{Learning} Sparse Support under Differential Privacy: Adaptive Algorithms and Minimax Limits''},
  author = {Jia Gu and T. Tony Cai},
  year   = {2026}
}

@article{Chang-JASA-2026,
  author    = {Jinyuan Chang and Lin Yang and Mengyue Zha and Wen-Xin Zhou},
  title     = {Adapting to Noise Tails in Private Linear Regression},
  journal   = {Journal of the American Statistical Association},
  pages     = {1--14},
  year      = {2026},
  publisher = {Taylor \& Francis}
}

@article{Cai-JRSSB-2026,
    author = {Cai, Zhanrui and Fan, Yingying and Gao, Lan},
    title = {Knockoff inference under privacy constraints},
    journal = {Journal of the Royal Statistical Society Series B: Statistical Methodology},
    pages = {qkag120},
    year = {2026},
    month = {08},
    issn = {1369-7412},
    eprint = {https://academic.oup.com/jrsssb/advance-article-pdf/doi/10.1093/jrsssb/qkag120/70506907/qkag120.pdf},
}

@misc{Barber-Duchi-2014,
  author        = {Barber, Rina Foygel and Duchi, John C.},
  title         = {Privacy and Statistical Risk: Formalisms and Minimax Bounds},
  year          = {2014},
  eprint        = {1412.4451},
  archivePrefix = {arXiv},
  primaryClass  = {stat.ML},
  note          = {arXiv:1412.4451}
}

@misc{cai2026cost,
  author        = {Cai, T. Tony and Wang, Yichen and Zhang, Linjun},
  title         = {The Cost of Privacy: Rates of Convergence for
                   Parameter Estimation with Differential Privacy},
  year          = {2026},
  eprint        = {1902.04495},
  archivePrefix = {arXiv},
  primaryClass  = {stat.ML},
  note          = {arXiv:1902.04495v6, corrected version}
}

@misc{Wang-Shen-arxiv-2026,
  title         = {{NetPTR}: Optimal Differentially Private Spectral Community Detection on Sparse Networks},
  author        = {Wanjie Wang and Tao Shen},
  year          = {2026},
  eprint        = {2606.26145},
  archivePrefix = {arXiv},
  primaryClass  = {cs.SI},
  note          = {arXiv:2606.26145},
  url           = {https://arxiv.org/abs/2606.26145}
}
\end{document}